\RequirePackage{fix-cm}
\documentclass[times]{article}          % twocolumn
\usepackage{graphicx}
\usepackage{amsmath}

\usepackage{stmaryrd}

\usepackage{caption}
\usepackage{subcaption}
\usepackage[colorlinks,citecolor=blue,linkcolor=blue,urlcolor=blue]{hyperref}

\usepackage{doi}
\usepackage[authoryear]{natbib}
\setcitestyle{authoryear,open={(},close={)}}

\usepackage{float}

\title{Minimum Feature Size Control in Level Set Topology Optimization via Density Fields}

\author{Jorge L. Barrera$^{1*}$, Markus J. Geiss$^2$, Kurt Maute $^3$ \\[12pt]
$^1$Computational Engineering Division, Lawrence Livermore National Laboratory,\\ 
7000 East Ave Livermore, 94550, CA, USA\\
$^2$OHB System AG, Manfred-Fuchs-Strasse 1, 82234 Wessling, Germany \\ 
$^3$Ann and H. J. Smead Department of Aerospace Engineering Sciences,\\ 
University of Colorado at Boulder, 3775 Discovery Dr, Boulder, CO 80303, USA\\
$^*$Corresponding author: barreracruz1@llnl.edu}

\begin{document}

%\date{Received: date / Accepted: date} % The correct dates will be entered by the editor

\maketitle

\begin{abstract}
A level set topology optimization approach that uses an auxiliary density field to nucleate holes during the optimization process and achieves minimum feature size control in optimized designs is explored. The level set field determines the solid-void interface, and the density field describes the distribution of a fictitious porous material using the solid isotropic material with penalization. These fields are governed by two sets of independent optimization variables which are initially coupled using a penalty for hole nucleation. The strength of the density field penalization and projection are gradually increased through the optimization process to promote a 0-1 density distribution. This treatment of the density field combined with a second penalty that regulates the evolution of the density field in the void phase, mitigate the appearance of small design features. The minimum feature size of optimized designs is controlled by the radius of the linear filter applied to the density optimization variables. The structural response is predicted by the extended finite element method, the sensitivities by the adjoint method, and the optimization variables are updated by a gradient-based optimization algorithm. Numerical examples investigate the robustness of this approach with respect to algorithmic parameters and mesh refinement. The results show the applicability of the combined density level set topology optimization approach for both optimal hole nucleation and for minimum feature size control in 2D and 3D. This comes, however, at the cost of a more advanced problem formulation and additional computational cost due to an increased number of optimization variables.
\end{abstract}

% -------------------------------------------------------------------
% -------------------------------------------------------------------
\section{Introduction} \label{sec:intro}
% -------------------------------------------------------------------
% -------------------------------------------------------------------

In topology optimization (TO), mesh-independent and manufacturable designs are achieved by controlling the length scale of structural members. Density-based TO methods possess robust and reliable approaches for this; see for example \cite{Liu2016} or \cite{Lazarov2016}. In contrast, less mature strategies are available for level set (LS) methods. In this paper, a novel approach that controls the minimum feature size in LS TO using density fields is developed and studied for linear elastic structural problems in two and three dimensions.

In density-based TO approaches, control of the minimum feature size is achieved by filtering and projecting the optimization variables; see, for example, \cite{Guest2009,Guest2004,Guest2009a,Sigmund2007,Sigmund2009}. The filtering scheme widens the area of influence of each design variable while projecting the densities alleviates the blurriness introduced by filtering. This strategy reliably achieves features with the desired length scale but allows for the existence of point hinges in optimal designs. Hence, alternatives such as robust optimization were introduced, see, for example, \cite{Schevenels2011,Asadpoure2011,fernandez2020imposing,Sigmund2011}. In this strategy, the filter size is perturbed to alter the shape of design by the desired minimum feature size, the worst-case configuration identified, and its performance optimized. This results in a nominal design that is insensitive to the imposed shape variations which is typically achieved by the minimum feature of the optimized nominal design being larger than desired minimum feature size.
%@ProfNew: I don't get this sentences. Insensitivity of the nominal design wrt feature size as in it doesn't care about it?
%@Jorge: Is this better? Still a bit convoluted for me, but it is what it is. 

Implicit LS methods solve some form of the Hamilton-Jacobi equation to track the evolution of external and internal boundaries in the optimization process (\cite{van2013level,sigmund2013topology}). For such LS methods, the skeleton approach has been proposed to control feature size (\cite{Guo2014, Xia2015, Liu2018}). Here, the skeleton is defined as the location of centers of maximum circles (in 2D) or spheres (in 3D) contained within the body that is defined by the zero LS isocontour (\cite{Montanari1968}). It is extracted from a signed distance field (SDF) constructed through reinitialization of the LSF when solving the Hamilton-Jacobi equation. Hence, the skeleton approach assumes a unit gradient of the level set field (LSF) to measure the feature size using the LS values of skeleton points (\cite{Alexandrov2005,Kimmel1993}). Despite commonly used, this approach suffers from convergence issues. Small perturbations in shape may lead to a drastically different skeleton, which prevents a smooth design evolution; see \cite{Allaire2016}.

In explicit LS methods, the parameters of the discretized LSF are defined as explicit functions of the optimization variables and updated by optimization algorithms. Recently, \cite{Geiss2019a} studied the skeleton approach using an explicit LS method. In his work, the skeleton is constructed from the Laplacian field of a SDF generated using the heat method of \cite{crane2017heat}. Control on the feature size is achieved following the work of \cite{Guo2014}. It was found that this approach requires a fine mesh to sufficiently resolve the Laplacian of the SDF, which leads to an increased computational cost. In addition, the SDF needs to be filtered to obtain sufficiently smooth higher-order spatial derivatives. This, in turn, leads to smearing of the LS skeleton which affects the accuracy of enforcing feature size. 

The work of \cite{Villanueva2016} and \cite{Coffin2016} proposed an alternative to the skeleton approach for explicit LS methods. Their strategy focuses on using explicit geometric distance measures to control the size of minimum features. However, their approach lacks generality as it is limited to controlling features of the size of an element of the mesh.

In this paper, the classical solid isotropic material with penalization (SIMP) approach is included into an explicit LS TO framework to achieve minimum feature size control. Similar ideas have been explored by various authors. 
In \cite{jansen2019explicit}, feature size control is achieved via geometric constraints on a density field, which is an extension of the approach presented by \cite{Dunning2018}. Furthermore, a combined LS-density approach was presented in \cite{geiss2018topology,geiss2019combined} where the LS and density fields are coupled while optimizing the material distribution within the solid phase. While this approach enables solving of multi-material TO problems, it provides no means of feature size control. Recently, \cite{JANSEN2020} introduced a hybrid approach for multi-scale TO using density and LS methods. In their work, the density field is used to represent homogenized material properties at a micro-scale while the LSF is used to describe the geometry at a macro scale. Minimum length scale control has also been explored by introducing two additional geometric constraints, as initially proposed by \cite{Zhou2015} for three-field SIMP problems. Most recently, \cite{Andreasen2020} introduced a hybrid LS-density TO scheme for 2D problems which provides length-scale control via a density approach and then performs shape optimization via the LS method. In their work, a density field is used as the design variable on which linear filtering is applied and ultimately a LSF is computed. While minimum feature size is achieved via the robust density approach (\cite{Sigmund2009}), holes are created via a separate heuristic nucleation scheme which requires additional parameter tuning. 
%<Markus: this statement is not clear; be more explicit here, specifically wrt to the work of Nils> 
%Due to only a single set of optimization variables, natural seeding of holes is not possible and a large dependency of the final designs on the initial guess is observed.
%@Kurt: I hope it's more clear now how we differ from Niels' work and others.

The TO approach proposed in this paper nucleates holes informed by a density field and enables minimum feature size control via the same density field. This approach is an extension of the work presented by \cite{Barrera2019Hole} addressing the issue of minimum feature size control in explicit LS TO. In the solid-void problems considered in this study, a LSF is used to distinguish between the two material phases through a crisp, well-defined interface. This LSF is parameterized by local shape functions and the nodal LS values are defined as explicit functions of LS optimization variables. The density field interpolates the material properties within the material phase using the SIMP scheme (\cite{bendsoe2004topology}) and is discretized on the same mesh as the LSF. This density field is obtained via linear filtering and projection of nodal density design variables, represented by a dedicated set of density optimization variables. The density field is used to nucleate holes in the LSF and to provide minimum feature size control on the explicit LS TO problems. The weak form of the governing equations is discretized by the extended finite element method (XFEM) \cite{belytschko2009review}. The proposed approach is studied for compliance and mass minimization problems in 2D and 3D, considering linear elasticity. The combined LS-density TO approach successfully demonstrates the ability to obtain feature size control and to yield mesh-independent designs. To achieve this, a penalization term for removing high density areas in the void phase needs to be employed and a globally regularized LSF is indispensable.

The remainder of this paper is organized as follows: Section \ref{sec:LSTO} introduces the LS and density problems in the proposed TO approach. Section \ref{sec:holeSeedAndFSCSec} details the strategies used for hole seeding and feature size control using the density method. The optimization framework and problem formulation are introduced in Section \ref{sec:LsOptFramework}. Numerical examples illustrating the properties of the proposed minimum feature size control approach are studied in Section \ref{sec:NumEx}, and conclusions and topics for future work are discussed in Section \ref{sec:Concl}.

% -------------------------------------------------------------------
% -------------------------------------------------------------------
\section{Topology optimization approach} \label{sec:LSTO}
% -------------------------------------------------------------------
% -------------------------------------------------------------------

\subsection{Explicit level set topology optimization} \label{subsec:expLsTopOpt}
% -------------------------------------------------------------------
Considering a two-phase configuration, a LSF defined on the design domain, $\Omega_D$, describes the geometry as follows:
% Equation
\begin{equation}\label{eq:LSDescription}
\begin{aligned}
\phi(\boldsymbol X)
\begin{cases}
> 0, ~~~~~& \forall~\boldsymbol X \in \Omega_I, \\
< 0, ~~~~~& \forall~\boldsymbol X \in \Omega_{II}, \\
= 0, ~~~~~& \forall~\boldsymbol X \in \Gamma_{I,II}.
\end{cases}
\end{aligned}
\end{equation}
The material domains of phases $I$ and $II$ are identified by $\Omega_I$ and $\Omega_{II}$, respectively, such that $\Omega_D = \Omega_I \cup \Omega_{II}$; and the zero LS isocontour, $\phi(\mathbf{X}) = 0$, defines the interface, $\Gamma_{I,II}$. In the solid-void problems studied in this paper, solid is assigned to phase I and void to phase II.

The LSF is discretized by local shape function and defined as an explicit function of a vector of LS optimization variables, 
${\boldsymbol s}^\phi := \{ {\boldsymbol s}^\phi \in \rm I\!R^{N_s} |~ \phi_{low} \leq {s}^\phi_i \leq \phi_{up}, i=1,...,N_s \}$.
Here, a LS optimization variable is assigned to each node of a structured mesh; thus, $N_s$ equals the number of nodes in this mesh. The optimization variables are updated via a nonlinear programming method. Note that this framework does not require solving a Hamilton-Jacobi-type equation, which is commonly needed in implicit LS optimization approaches; see \cite{van2013level} for details.

The LS optimization variables, ${\boldsymbol{s}}^\phi$, are filtered following the formulation presented in \cite{kreissl2012levelset} to generate a vector of filtered LS coefficients, 
$\boldsymbol{\hat{s}}^\phi := \{ \boldsymbol{\hat{s}}^\phi \in \rm I\!R^{N_s} |~\phi_{low} \leq \hat{s}^\phi_i \leq \phi_{up}, i=1,...,N_s \}$.
The purpose of this distance-based linear filter is to enhance the smoothness and the convergence of the design problem. A filtered LS coefficient, $\hat{s}^\phi_i$, at node $i$, is defined as a function of its neighboring LS optimization variables, ${s}^\phi_j$, at nodes $j$, as follows:
% Equation
\begin{equation}\label{eq:LSLinearFilter1}
\begin{aligned}
\hat{s}_i^\phi = {F_{ij}(r_f^\phi)} \ s_j^{\phi}
\end{aligned}
\end{equation}
where the operator ${F}_{ij}$ is function of the LS filter radius, $r_f^\phi$ as follows:
% Equation
\begin{equation}\label{eq:LSLinearFilter2}
\begin{aligned}
{F}_{ij}(r_f^\phi) = w_{ij}(r_f^\phi) / \sum^{N_{r^\phi_f}}_j{w_{ij} (r_f^\phi)}.
\end{aligned}
\end{equation}
The parameter $N_{r^\phi_f}$ denotes the number of nodes within $r^\phi_f$, and the weights $w_{ij}$ are defined as: 
% Equation
\begin{equation}\label{eq:LSLinearFilter3}
\begin{aligned}
w_{ij}(r_f^\phi) = \max(0,r^\phi_f-|\mathbf{X}_i-\mathbf{X}_j|),
\end{aligned}
\end{equation}
where the expression $|\mathbf{X}_i-\mathbf{X}_j|$ represents the Euclidean distance between nodes i and j.
The LSF is linearly interpolated by the filtered LS coefficients, $\boldsymbol{\hat{s}}^\phi$, on linear quadrilateral and hexahedral meshes in 2D and 3D, respectively, using:
% Equation
\begin{equation}\label{eq:discrLsField}
\begin{aligned}
\phi( \mathbf{X})
= 
\displaystyle\sum_{i=1}^{N^e} \mathcal{N}_i ( \mathbf{X}) ~ \hat{s}_i^\phi.
\end{aligned}
\end{equation}
Linear shape functions are denoted by $\mathcal{N}( \mathbf{X})$, and the LSF is defined as
$\phi ( \mathbf{X}):=\{ \phi ( \mathbf{X}) \in H^1(\Omega_D) |~ \phi_{low} \leq a^\phi \leq \phi_{up}:  a^\phi \in \mathcal{\hat{S}^\phi} ( \mathbf{X})  \}$, 
with $H^1$ denoting the Sobolev space.
The upper and lower bounds of the LSF are denoted by the parameters $\phi_{low}$ and $\phi_{up}$, respectively.

\subsection{Density method} \label{subsec:densMeth}
% -------------------------------------------------------------------
In addition to the $LS$ optimization variables introduced above, we define a second vector, 
${\boldsymbol s}^\rho := \{ {\boldsymbol s}^\rho \in \rm I\!R^{N_s} |~ 0 \leq {s}^\rho_i \leq 1, i=1,...,{N_s} \}$, of $density$ optimization variables to control a fictitious density field as described in this Section. Hence, the proposed framework considers a set composed of nodal LS and density optimization variables, i.e., ${\boldsymbol s} =[ {\boldsymbol s}^\phi, {\boldsymbol s}^\rho ]$.

Similar to the filtering described in the previous Section for the LS optimization variables, here a vector of filtered density coefficients,
${\boldsymbol{\hat{s}}}^\rho := \{ {\boldsymbol{\hat{s}}}^\rho \in \rm I\!R^{N_s} |~ 0 \leq \hat{s}^\rho_i \leq 1, i=1,...,{N_s} \}$,
 is generated using the following: 
% Equation
\begin{equation}\label{eq:DensLinearFilter}
\begin{aligned}
\hat{s}_i^\rho = {F_{ij}(r_f^\rho)} \ s_j^{\rho},
\end{aligned}
\end{equation}	
with $F_{ij}$ defined in Eq.~\ref{eq:LSLinearFilter2}.
Note that the density filter radius, $r^\rho_f$ in this expression is not required to be the same as the LS filter radius, $r_f^\phi$ in Eq.~\ref{eq:LSLinearFilter1}. The set $\boldsymbol{\hat{s}}^\rho$ is used to linearly interpolate a filtered density design variable field,
$\mathcal{\hat{S}^\rho} ( \mathbf{X}):=\{ \mathcal{\hat{S}^\rho} ( \mathbf{X}) \in H^1(\Omega_D) |~ 0 \leq a^\rho \leq 1:  a^\rho \in \mathcal{\hat{S}^\rho} ( \mathbf{X})  \}$. This field is also linearly interpolated using:
% Equation
\begin{equation}\label{eq:discrLsField}
\begin{aligned}
\mathcal{\hat{S}^\rho} ( \mathbf{X}) 
= 
\displaystyle\sum_{i=1}^{N^e} \mathcal{N}_i ( \mathbf{X}) ~ \hat{s}_i^\rho.
\end{aligned}
\end{equation}

It has been shown by \cite{Guest2004,sigmund2007morphology,guest2009imposing,wang2011projection}, among others, that TO using density methods can provide certain control on the minimum feature size by using a large penalization in combination with projection of the density field. 
Thus, $\mathcal{\hat{S}^\rho}(\mathbf{X})$ is projected to obtain the nodal, filtered, projected fictitious density field, $\hat\rho(\mathbf{X})$, hereinafter referred to as projected density field for brevity. To this end, the threshold projection formulation of \cite{wang2011projection} is adopted:
% Equation:
\begin{equation}\label{eq:SIMPProj}
\begin{aligned}
	\hat\rho( \mathbf{X})
	= 
	\frac{
	\tanh(\gamma_{pr} (\mathcal{\hat{S}^\rho} ( \mathbf{X}) - \tau_{pr}) + \tanh(\gamma_{pr}\tau_{pr})
	}{
	\tanh(\gamma_{pr} (1 - \tau_{pr}) + \tanh(\gamma_{pr}\tau_{pr})
	}.
\end{aligned}
\end{equation} 
The parameters $\tau_{pr}$ and $\gamma_{pr}$ represent the projection threshold and sharpness, respectively. For increasing values of $\gamma_{pr}$, this formulation projects any value above $\tau_{pr}$ to 1, and the values below to 0. In this work, $\tau_{pr}$ is set to a value sufficiently close to zero, i.e. $\tau_{pr}=0.001$ to ensure minimum length scale on the solid phase, emulating the strategy presented in \cite{Guest2004}. Note that setting $\tau_{pr}$ to a higher value (e.g., $\tau_{pr}=0.5$) would impede minimum length scale control, as demonstrated in \cite{wang2011projection}. Also, methods such as the robust approach (\cite{schevenels2011robust}) are needed to mitigate the appearance of point hinges in converged designs. 

\subsubsection{Material interpolation} \label{subsubsec:matInterpSIMP} 
% -------------------------------------------------------------------
In the linear elastic problems studied in Section \ref{sec:NumEx}, material properties are interpolated using the SIMP approach. The material density, $\theta(\mathbf{X})$, is defined as a linear function of the projected density field as follows:
% Equation
\begin{equation}\label{eq:varyingMatDensForm}
\begin{aligned}
	\theta(\mathbf{X}) = \theta_0  ~ \hat\rho(\mathbf{X}).
\end{aligned}
\end{equation}
In addition, the Young's modulus, $E(\mathbf{X})$, is interpolated using the following power law:
% Equation
\begin{equation}\label{eq:varyingMatPropForm}
\begin{aligned}
	E(\mathbf{X}) = E_0 \left(\hat\rho(\mathbf{X}) \right)^{\beta_{\rho}}.	
\end{aligned}
\end{equation}
The parameters $\theta_0$ and $E_0$ denote the bulk material density and Young's modulus, respectively. The $\beta_{\rho}$ parameter symbolizes the SIMP exponent. The nodal densities are averaged within each element, and an element-wise constant density is used for interpolation purposes.
%%% Are you also discussing somewhere how this is done for intersected XFEM elements? I think it would be worth wile to state how this is interpolated to the centroid and then assigned constant for each tri/tet
%<Jorge: Yes, I agree we should discuss this. However to this we would have to cover in more detail XFEM >
%<Prof: Should I add an xfem section so the paper is self-contained then? No need to add anything
Note that, despite the projected density field is defined in the entire design domain, the physical response is computed only in the solid phase, $\Omega_{I}$. 

\subsubsection{Continuation scheme for density method} \label{subsubsec:SIMPCont} 
% -------------------------------------------------------------------
Both the projection of $\hat\rho(\mathbf{X})$ defined in Eq.~\ref{eq:SIMPProj} and the SIMP exponent in Eq.~\ref{eq:varyingMatPropForm} are updated throughout the course of the optimization process using a continuation scheme. 
The projection sharpness, $\gamma_{pr}$, is set initially to a low value, i.e. $\gamma_{pr}=\gamma_{pr}^{0}=0.001$, to eliminate the effect of the projection at the beginning of the optimization process; and is gradually increased to $\gamma_{pr}=\gamma_{pr}^{f}=40.0$ to strengthen its effect and promote a 0-1 density field. The expression used reads:
% Equation:
\begin{equation}\label{eq:projContFunc}  
\begin{aligned}
\gamma_{pr} 
= 
	\begin{cases}
	\displaystyle \gamma_{pr}^{0} 
	+ 
	(\gamma_{pr}^{f}-\gamma_{pr}^{0})\left( \frac{\mathcal{D}_{it}}{\mathcal{D}_{c}}\right)^{\eta_{\gamma_{pr}}}, & \forall~\mathcal{D}_{it} \leq \mathcal{D}_{c} \\
	\gamma_{pr}^{f}, & otherwise,
	\end{cases}	
\end{aligned}
\end{equation}
where $\mathcal{D}_{it}$ is the design iteration index. The parameter $\mathcal{D}_{c}$ is the number of design iterations for which the continuation is active. The continuation step size, $\mathcal{D}_{st}$, is assumed to be constant. The projection is kept at its maximum value until the optimization problem converges (i.e., for $\mathcal{D}_{it} > \mathcal{D}_{c}$). The increase of $\gamma_{pr}$ between continuation steps is controlled by the exponent $\eta_{\gamma_{pr}}$, which is set to $\eta_{\gamma_{pr}}=2.0$. Similarly, the SIMP exponent $\beta_{\rho}$ is updated using the following expression:
% Equation:
\begin{equation}\label{eq:simpExpContFunc} 
\begin{aligned}
\beta_{\rho}
= 
	\begin{cases}
	\displaystyle 
	\beta_{\rho}^{0} 
	+ 
	(\beta_{\rho}^{f}-\beta_{\rho}^{0})\left( \frac{\mathcal{D}_{it}}{\mathcal{D}_{c}}\right)^{\eta_{\beta_{\rho}}}, & \forall~\mathcal{D}_{it} \leq \mathcal{D}_{c} \\
	\beta_{\rho}^{f}, & otherwise.
	\end{cases}	
\end{aligned}
\end{equation}
The parameters $\beta_{\rho}^{0}$ and $\beta_{\rho}^{f}$ are set to 2.0 and 12.0, respectively, and the exponent $\eta_{\beta_{\rho}}$ is also set to 2.0. The values for $\gamma_{pr}^{0}$, $\gamma_{pr}^{f}$, $\beta_{\rho}^{0}$, and $\beta_{\rho}^{f}$ were determined through numerical experiments. These values are sufficient to achieve both a high projection parameter (to counteract the smoothing introduced by the linear density filter) and high SIMP exponent (to promote convergence to a 0-1 design) without excessively increasing the nonlinearity of the problem. Figure~\ref{fig:continuationScheme_forDensMeth} 
illustrates the continuation strategy used for updating both $\gamma_{pr}$ and $\beta_{\rho}$ (i.e., Eqs.~\ref{eq:projContFunc} and \ref{eq:simpExpContFunc} , respectively). Note that an exponent greater than 1.0 attenuates changes in the updated parameters at early stages of the optimization process.

% Figure
\begin{figure}[h!]
	\centering
	\includegraphics[width=0.5\linewidth]{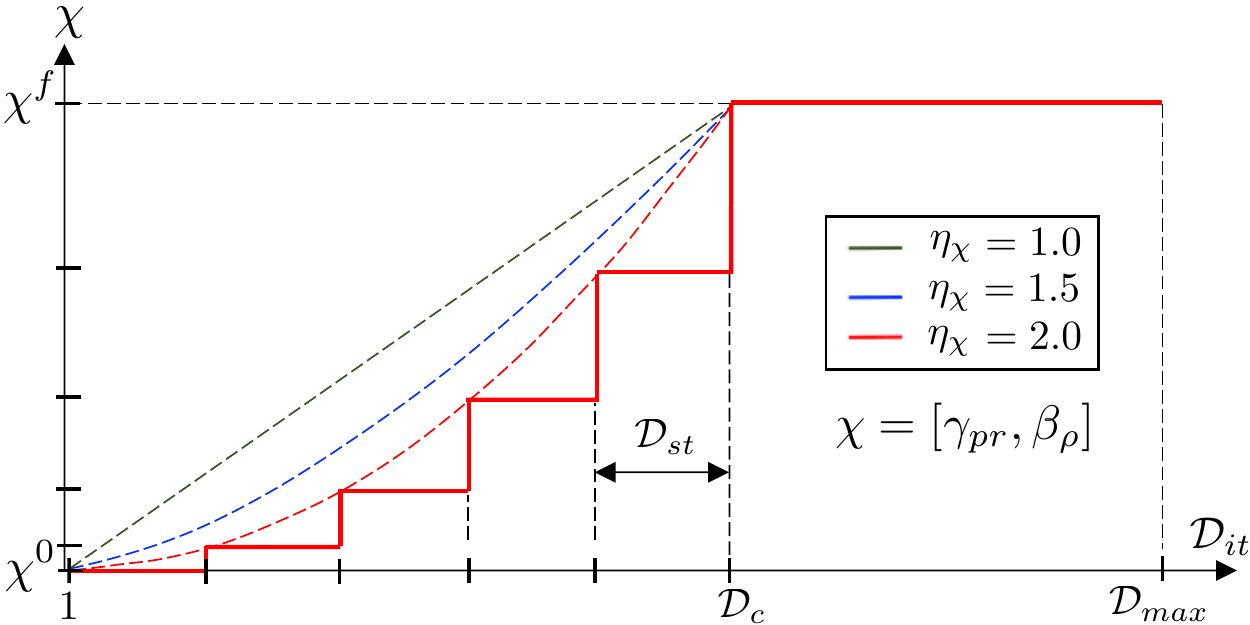}
	\caption{Continuation strategy for projection sharpness and SIMP exponent.}
	\label{fig:continuationScheme_forDensMeth} 
\end{figure}

% -------------------------------------------------------------------
% -------------------------------------------------------------------
\section{Hole nucleation and minimum feature size control via density methods} \label{sec:holeSeedAndFSCSec}
% -------------------------------------------------------------------
% -------------------------------------------------------------------

This section introduces an approach that aims to alleviate two shortcomings of the explicit LS optimization framework outlined above: (i) dependency of the optimization result on the initial hole seeding, and (ii) lack of control of the minimum length scale. In this approach, a density field that is independent of the design LSF is constructed. The mechanism that affords hole nucleation during the optimization process uses the hole seeding strategy summarized in Section~\ref{subsec:holeSeeding}. 
The strategy for minimum feature size control is discussed in Section~\ref{subsec:fscStrategy}. 

\subsection{Hole seeding strategy} \label{subsec:holeSeeding}
% -------------------------------------------------------------------
Building on the work of \cite{Geiss2019a,Barrera2019,Barrera2019Hole}, hole nucleation is established by coupling the LS and density fields such that: (i) high densities are present only in the solid phase, where the LSF is positive; and (ii) densities close to zero exist only in the void phase, where the LSF is negative. This is achieved by introducing a differentiable penalty term that couples the LSF to the projected density field into the objective function.

\subsubsection{Hole seeding penalty} \label{subsec:holeSeedingPen}
% -------------------------------------------------------------------
The smooth penalty formulation that relates the LSF, $\phi(\mathbf{X})$, and the projected density field
, $\hat\rho(\mathbf{X})$, reads:
% Equation
\begin{equation}\label{eq:TFCFomAllDomain}
\begin{aligned}
	{p}_{\hat\rho\phi}(\mathbf{X}) &=
	\begin{cases}
	\displaystyle\frac{
	\left\{ \left[ \max \left( 0, \displaystyle\frac{\phi -\phi_{th_{hs}}}{\phi_{up}-\phi_{th_{hs}} } \right) \right]^2 
	+ \xi^2 \right\}^{\frac{1}{2}} - \xi
	}{
	\left( 1 + \xi^2 \right)^{\frac{1}{2}}  - \xi
	}, 			
	\\ 
	~~~~~~~~~~~~~~~~~~~~~~~~~~~~~~~~~~~~~~\forall~\hat\rho<\rho_{th_{hs}}, \\
	0, ~~~~~~~~~~~~~~~~~~~~~~~~~~~~~~~~~~~otherwise;	
	\end{cases} \\
\end{aligned}
\end{equation}
where the parameter $\xi$ smooths the transition of the penalty and is set to $\xi=0.1$. Fig.~\ref{fig:TFCScheme} shows a schematic of this coupling penalty formulation. The LS threshold, $\phi_{th_{hs}}$, is set to a value below zero, i.e. $\phi_{th_{hs}}=0.10 \phi_{low}$ with $\phi_{low}<0$, to ensure that the penalty is active until after the LSF crosses the zero isocontour. The hole seeding density threshold, $\rho_{th_{hs}}$, increases from an initial low value, $\rho^0_{th_{hs}}$, to a value close to 1.0, i.e. $\rho^{f}_{th_{hs}}=0.98$, during the optimization process.

In the original strategy described in \cite{Barrera2019Hole}, the $\rho_{th_{hs}}$ threshold is decreased to zero. In contrast, here a gradual increase of this threshold is preferred since it widens the area of influence of the coupling penalty, and thus promotes the nucleation of holes in regions dominated by intermediate densities at a later stage in the design process. The extended hole nucleation effect together with projecting the density field largely mitigates the presence of intermediate densities in the solid phase away from the interface. 

Another difference with respect to the work of \cite{Barrera2019Hole} consists of the treatment of the density field used to interpolate material properties. Previously, the density field was gradually shifted through the optimization process until a uniform field of $1.0$ was obtained in the entire design domain. The goal of the shifting scheme was to guarantee that the optimization problem transitioned from a pure density problem to a pure LS problem. That is no longer the case in the approach presented here. The design is sensitive to both the density and LS components during the entire optimization process. A discussion of these and other differences between current and previous work are provided in Section \ref{subsubsec:discussion}. 
% Figure
\begin{figure}[t!]
	\centering
	\includegraphics[width=0.5\linewidth]{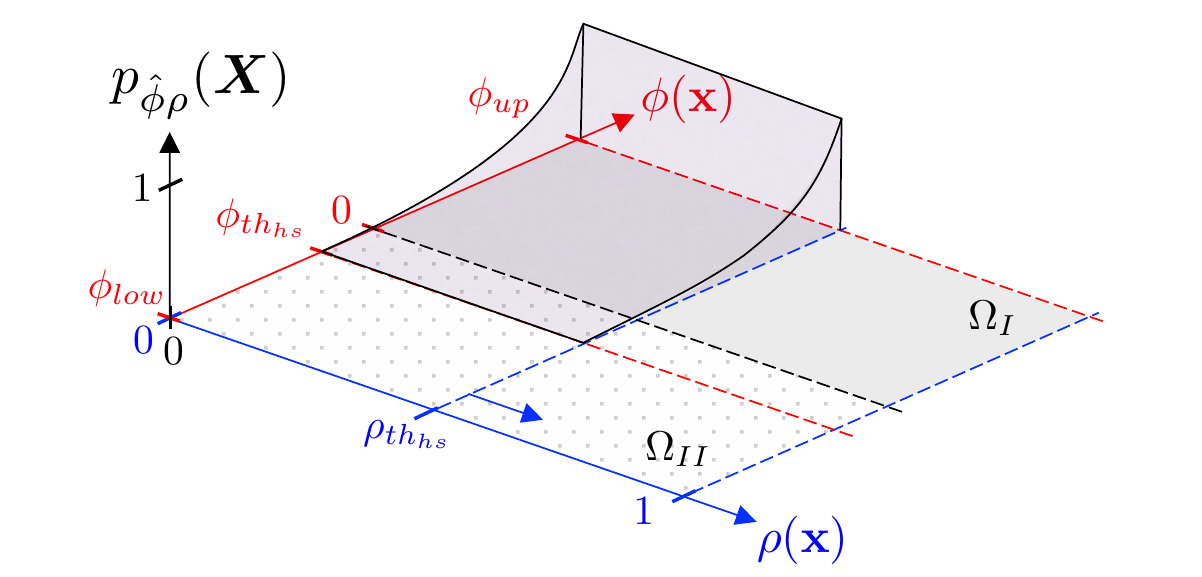}
	\caption{Hole seeding coupling penalization plotted over density and level variables.}
	\label{fig:TFCScheme} 
\end{figure}

\subsubsection{Continuation scheme for hole seeding strategy} \label{subsubsec:contRhoTh} 
% -------------------------------------------------------------------
A continuation scheme similar to the one presented in Section \ref{subsubsec:SIMPCont} is used to update the $\rho_{th_{hs}}$ threshold. This parameter is increased every $\mathcal{D}_{st}$ design iterations using:
% Equation
\begin{equation}\label{eq:expFuncContDec}
\begin{aligned}
	\rho_{th_{hs}} 
	&=
	\begin{cases}
	\displaystyle  \rho_{th_{hs}}^0 + \Delta\rho \left( \frac{\mathcal{D}_{it}}{\mathcal{D}_{c}}\right)^{\eta_{\rho_{th_{hs}}}},  ~~~~~~~~ \forall~\mathcal{D}_{it} \leq \mathcal{D}_{c}, \\
	\rho_{th_{hs}}^{f},   ~~~~~~~~~~~~~~~~~~ \forall~\mathcal{D}_{c} \leq \mathcal{D}_{it}  \leq \mathcal{D}_{c} + \mathcal{D}_{st}, \\
	0.0,  ~~~~~~~~~~~~~~~~~~~~~~~~~~~~~~~~~~~~~~otherwise,	
	\end{cases}
\end{aligned}
\end{equation}
with
% Equation
\begin{equation}\label{eq:expFuncContDecDeltaRho}
\begin{aligned}
   \Delta\rho &= \rho_{th_{hs}}^{f}- \rho_{th_{hs}}^0,
\end{aligned}
\end{equation}
as depicted in Fig.~\ref{fig:continuationScheme}. The parameter $\rho_{th_{hs}}^0$ is chosen to be a fraction of the initial prescribed homogenous density in the design domain. The exponent $\eta_{\rho_{th_{hs}}}$ is set to 2.0. The $\rho_{th_{hs}}$ threshold is set to zero for $\mathcal{D}_{it}>\mathcal{D}_{s}$ ($\mathcal{D}_{s}=\mathcal{D}_{c} + \mathcal{D}_{st}$) to deactivate the hole seeding penalty at a later stage of the optimization process and use the projected density field only to provide control on the minimum feature size allowed.
 % Figure
\begin{figure}[h!]
	\centering
	\includegraphics[width=0.5\linewidth]{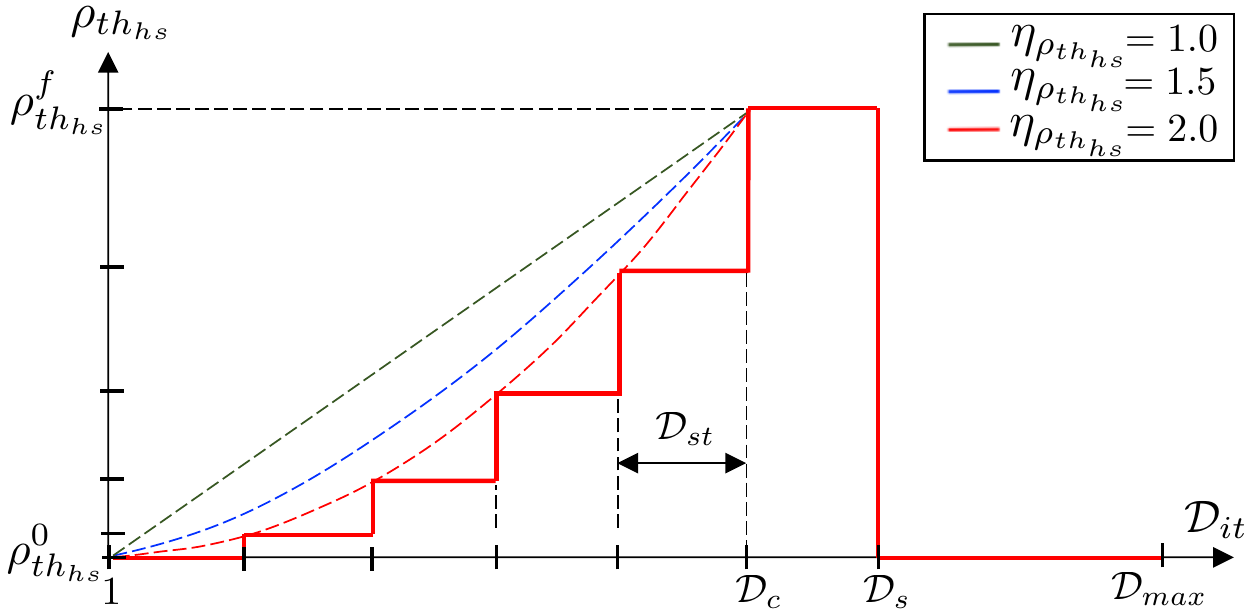}
	\caption{Continuation strategy for hole seeding density threshold.}
	\label{fig:continuationScheme} 
\end{figure}

\subsection{Feature size control strategy} \label{subsec:fscStrategy} 
\begin{figure*}[ht!]
	\centering
	\includegraphics[width=0.95\linewidth]{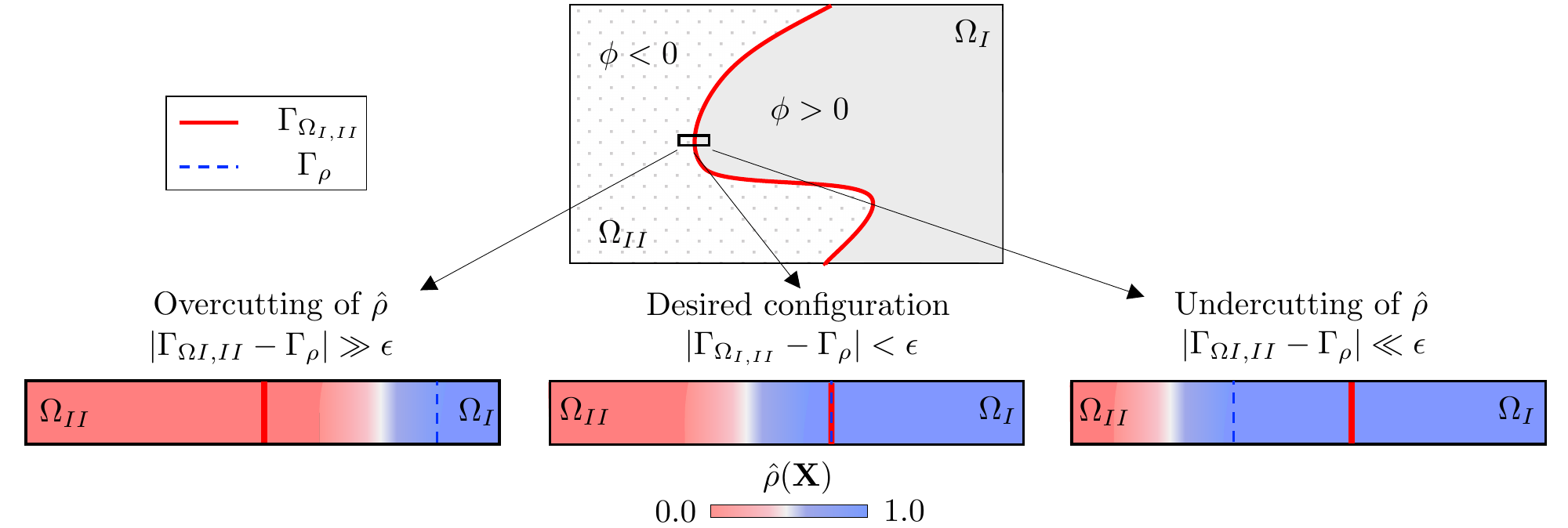}
	\caption{Insets with zoomed-in regions of vicinity of interface for three configurations: (lefft) overcutting, (center) desired cut and (right) undercutting of the density field by the interface.}
	\label{fig:DensityIssueInVoidRegion} 
\end{figure*}

% -------------------------------------------------------------------
The explicit LS framework introduced here enforces the minimum length scale provided by a density method.  To achieve this, the configuration 
illustrated at the inset in the middle of Fig.~\ref{fig:DensityIssueInVoidRegion} is considered. The $\Gamma_{\rho}(\boldsymbol X)$ curve/surface in 2D/3D identifies the material and void phases using the projected density field. Here, the zero isocontour of the LSF (i.e., the interface) tracks this contour within a tolerance $\epsilon$. The $\Gamma_{\rho}(\boldsymbol X)$ reference contour is constructed such that it separates regions of high densities $(\approx1.0)$ from regions of intermediate and low densities.  This contour is defined only for explanation purposes, but it is not explicitly used in the proposed framework.

% Figure

Attaining the scenario mentioned above while having two independent fields is far from trivial. As is demonstrated later in Section \ref{sec:NumEx}, setting algorithmic parameters inappropriately can result in configurations where the interface is away from $\Gamma_{\rho}(\boldsymbol X)$, i.e. $|\Gamma_{\Omega_{I,II}} - \Gamma_\rho| \gg \epsilon$. In such cases, the interface would either ``overcut" or ``undercut" the projected density field, as shown at the left and right insets of Fig.~\ref{fig:DensityIssueInVoidRegion}.  These two configurations suggest a disconnection between the LSF and projected density field, which needs to be avoided. Otherwise, the optimization framework is unable to reliably provide control of the minimum feature size. Overcutting implies that regions of low density are located in the solid phase at convergence. Extreme overcutting, i.e., large regions of low densities in the solid phase are unlikely to occur because both the density and LS problems are governed by the same optimization formulation and will remove material similarly. However, the interface may converge in the transition region to benefit from weak material near the interface. On the other hand, undercutting, i.e., regions of high densities in the void phase, can be mitigated by forcing the projected density field to be zero in the void phase. A penalty to achieve this is proposed next.

\subsubsection{Void domain density removal penalty} \label{subsec:VDDRPen}
%@Jorge: Please, consider using a different expression than \93clean up\94; this does not sound scientific; if you decide to change the expression please, fix it here and throughout the text.
% -----------------------------------------------------------------------------
In conventional SIMP and SIMP-like approaches, the material response is evaluated in a fixed design domain throughout the optimization process. In contrast, the approach described here restricts the effect of the density method to only the solid phase $\Omega_{I}$, which is defined by a moving interface. Hence, performance measures, i.e., criteria used to formulate objective and constraints, are insensitive to changes in the projected density field in the void phase. The evolution of this field is not determined by nor affects the optimization problem when $\phi(\mathbf{X})<0$. As a result, in absence of a mechanism to keep the densities in the void phase close to zero, the feature size control reference contour, $\Gamma_{\rho}(\boldsymbol X)$, is lost. To overcome this issue, an additional void domain density removal (VDDR) penalty, ${p}_{\hat\rho^{\Omega_{II}}}$, is introduced. This penalty is active only in the void phase and does neither affect the performance measures nor the physical response.
Undercutting of the density field is averted by the following formulation of the VDDR penalty:
% Equation
\begin{equation}\label{eq:CleaupPenForm} 
\begin{aligned}
	{p}_{\hat\rho^{\Omega_{II}}}(\mathbf{X}) &
	=
	\begin{cases}
	\mathcal{\hat{S}^\rho} ( \mathbf{X}), &
	\forall~\phi\leq\phi_{th_{fs}} 
	\\
	\displaystyle 
	0, & 			 
	otherwise.
	\end{cases} \\
\end{aligned}
\end{equation}
This formulation penalizes the filtered density field, $\mathcal{\hat{S}^\rho} (\mathbf{X})$. In our experience, a smoother and quicker removal of high densities is achieved if the VDDR penalty is applied on the filtered density field instead of the projected density field (i.e., $\hat\rho(\boldsymbol X)$). The threshold $\phi_{th_{fs}}$ is set to a value such that it limits the influence of the penalty to the void side with an offset (from the zero isocontour), i.e., $\phi_{th_{fs}}=[0.25-0.75]\phi_{low}$. Note that enforcing this penalty too close to the interface, within an element length, can result in overcutting of the density. This is demonstrated in Example \ref{subsec:2dBeamInterplay}.

% Figure
\begin{figure}[h!]
	\centering
	\includegraphics[width=0.5\linewidth]{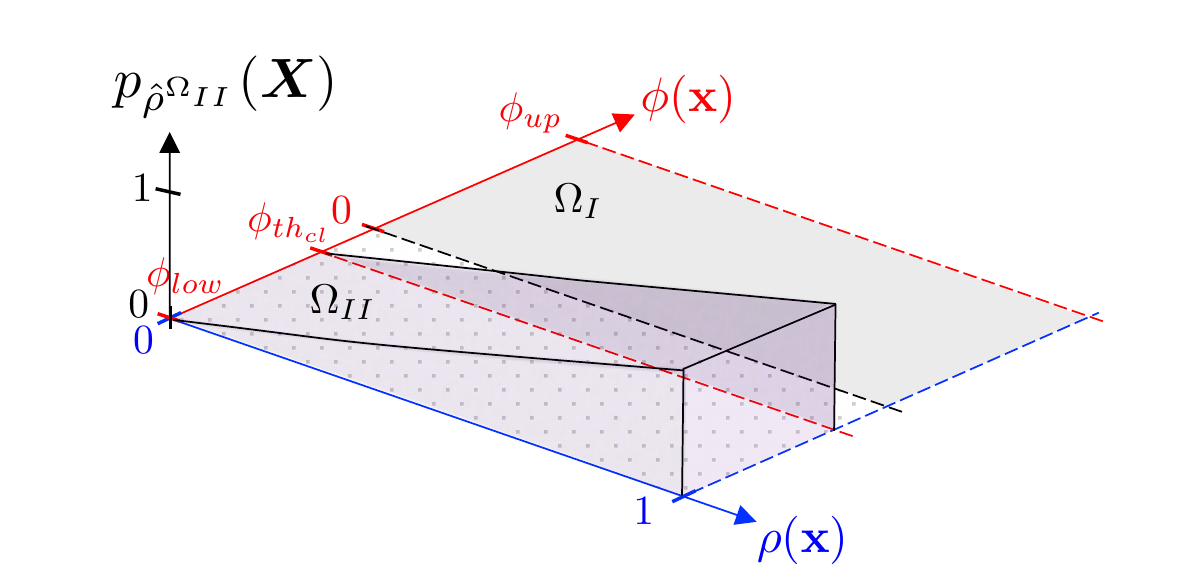}
	\caption{VDDR penalty in the void phase plotted over density and LS variables.}
	\label{fig:VDDRPen} 
\end{figure}

\subsection{Discussion} \label{subsubsec:discussion}
% -------------------------------------------------------------------
In the proposed TO framework, the design is driven by both the LS and density optimization variables until convergence. Thus, the optimization problem does not transitions from a pure density problem to a pure LS problem, as it was the case in previous work; see \cite{Barrera2019}. Conceptually, the framework described here resembles a combined LS-density approach rather than a LS approach that uses a density field to mitigate some of its disadvantages. Consequently, since the projected density field can be optimized within the solid domain, a uniform density field of 1.0 cannot be guaranteed in the optimized design. The hole seeding penalty alleviates this issue by promoting hole nucleation in regions with intermediate densities; see Section \ref{subsec:holeSeedingPen}. However, designs with intermediate densities are still possible if the optimization problem benefits from them (\cite{maute2017topology}). 

Unlike in previous works where the density field in the void phase had no influence on the optimization problem, here low densities are necessary in regions of the design domain where $\phi<0$ to prevent undercutting. Nevertheless, this can be easily achieved by a penalty term; see Section \ref{subsec:fscStrategy}.
Also, since the length scale is controlled only by the filter radius of the density optimization variables, $r^\rho_f$, the LS filter radius is not required to change as function of the desired minimum feature size.

A total of three parameters, i.e. $\phi_{th_{hs}}$, $\gamma_{pr}$, and $\beta_{\rho}$, are updated simultaneously through the optimization process every $\mathcal{D}_{st}$ design iterations; see Eqs.~\ref{eq:projContFunc}, \ref{eq:simpExpContFunc}, and \ref{eq:expFuncContDec}. As a result, the functionality provided by the projected density field transitions from first nucleating holes to finally controlling the minimum feature size. The $\phi_{th_{fs}}$ threshold of the VDDR penalty formulated in Section \ref{subsec:VDDRPen} determines an offset between the interface and the region where this penalty is active. A LSF with irregular gradients (i.e., locally too flat or steep) in the vicinity of the interface cannot be used to adequately enforce the VDDR penalty. Hence, a regularized LSF with a close to uniform gradient norm is critical. In this work, this is achieved by the LS regularization scheme presented in Section \ref{sec:LsOptFramework}. It is possible to avoid explicitly enforcing a VDDR penalty ${p}_{\hat\rho^{\Omega_{II}}}$ for a specific set of optimization problem formulations. This is explored in Section \ref{subsec:massMin3DBeam} with a mass minimization TO problem that defines the mass in the entire domain. In that setup, the void phase is part of the minimization problem, and thus high densities are penalized in the void phase by construction.

% -------------------------------------------------------------------
% -------------------------------------------------------------------
\section{Optimization problem formulation} \label{sec:LsOptFramework} 
% -------------------------------------------------------------------
% -------------------------------------------------------------------

The optimization problems considered in this work are formulated as follows: 
% Equation
\begin{equation}\label{eq:TFCOptProbSetup}
\begin{aligned}
\underset{s}{\min}~z(\mathbf{s}, \mathbf{u (\mathbf{s})})  
& = 
 w_1 ~ \mathcal{F}(\mathbf{s}, \mathbf{u(\mathbf{s})})
+ 
w_2  ~ P_{Per}(\mathbf{s})
\\
& 
+ 
w_3 ~ P_{Reg}(\mathbf{s})
+ 
w_4 ~ P_{\rho\phi}(\mathbf{s})
\\
& 
+ 
w_5 ~ P_{\hat\rho^{\Omega_{II}}}(\mathbf{s})
\\
s.t.: ~~~
&g_i(\mathbf{s}, \mathbf{u(\mathbf{s})}) \leq 0, i=1,...,N_g.
\end{aligned}
\end{equation}
The objective, $z$, is minimized over the vector of admissible optimization variables, ${\boldsymbol s} =[ {\boldsymbol s}^\phi, {\boldsymbol s}^\rho ]$. The state variables, $\mathbf{u}$, with $\mathbf{u} \in \rm I\!R^{N_{u}}$ $N_u$ being the number of state variables, are defined as functions of the optimization variables through the governing equations. The first component of the objective represents the performance measure to be minimized, $\mathcal{F}$ (e.g., strain energy, mass). The second term is the normalized perimeter control penalty, 
% Equation
\begin{equation}\label{eq:perPenForm}
\begin{aligned}
P_{Per}
= 
\frac{\displaystyle\int_{\Gamma_{I,II}} dA}{\displaystyle\int_{\Gamma_{D}} dA},
\end{aligned}
\end{equation}
which prevents the emergence of irregular geometric features. The perimeter of the design domain is denoted by $\Gamma_{D}$. The LS regularization penalty, $P_{Reg}$, is included in the objective to avoid spurious oscillations in the LSF. Following the regularization scheme of \cite{geiss2019regularization}, a truncated signed distance field used as target field, $\tilde\phi$, is constructed in the entire design domain. The target field is enforced through a penalty in the objective, $P_{Reg}$, of the following form:
% Equation
\begin{equation}\label{eq:LsRegEq}
\begin{aligned}
	P_{Reg} = 
	w_{\phi} \frac{\displaystyle\int_{\Omega_D} (\phi-\tilde\phi)^2 dV}{\displaystyle\int_{\Omega_D} \tilde\phi_{Bnd}^2 dV} + 
	w_{\nabla\phi} \frac{\displaystyle\int_{\Omega_D} | \nabla \phi- \nabla\tilde\phi|^2 dV}{\displaystyle\int_{\Omega_D} dV},
\end{aligned}
\end{equation}
where $\tilde\phi_{Bnd}$ is the difference between the upper, $\tilde\phi_{up}$, and lower, $\tilde\phi_{low}$, bounds in the target LSF. The weights $w_{\phi}$ and $w_{\nabla\phi}$ are kept constant and equal to one in the entire design domain. The fourth component denotes the normalized hole seeding coupling penalty,
% Equation
\begin{equation}\label{eq:CouplingPenForm}
\begin{aligned}
P_{\hat\rho\phi}
= 
\frac{\displaystyle\int_{\Omega_D^0} {p}_{\hat\rho\phi} dV}{\displaystyle\int_{\Gamma_{D}}  dA},
\end{aligned}
\end{equation}
with ${p}_{\hat\rho\phi}(\mathbf{X})$ as formulated in Eq.~\ref{eq:TFCFomAllDomain}. The last term is the normalized VDDR penalty to prevent undercutting, and is defined as
% Equation
\begin{equation}\label{eq:VDDRPenForm}
\begin{aligned}
P_{\hat\rho^{\Omega_{II}}}
= 
\frac{\displaystyle\int_{\Omega_D^0} {p}_{\hat\rho^{\Omega_{II}}}dV}{\displaystyle\int_{\Gamma_{D}}  dA}.
\end{aligned}
\end{equation}
The field ${p}_{\hat\rho^{\Omega_{II}}}(\mathbf{X})$ is computed using Eq.~\ref{eq:CleaupPenForm}.
The weights $w_1$, $w_2$, $w_3$, and $w_4$ in Eq.~\ref{eq:TFCOptProbSetup} are chosen such that all penalty contributions in the objective are significantly lower ($\approx 1-5\%$) than $\mathcal{F}$. Since a strong perimeter control penalty might prevent the nucleation of small holes (\cite{wang2007hole}), its contribution is kept below $1\%$ by gradually updating the $w_2$ weight, following the continuation scheme presented for the density method; see Section \ref{subsubsec:SIMPCont}. The weight of the VDDR penalty $w_5$ does not need special considerations since it affects only the density field in the void domain. Numerical examples have shown that using a value between $0.01$ and $10.0$ provides a relatively fast reduction of densities in the void phase.

The design needs to satisfy a set of $N_g$ problem dependent inequality constraints, $[g_1,...,g_{N_g}]$ (e.g., target mass, maximum allowable stress, maximum eigenvalue). The constraints are defined for each problem studied in Section \ref{sec:NumEx}. Note that the integral form of the coupling and VDDR penalties can be alternatively formulated via agglomerated nodal constraints such that they are enforced per node of the discretized domain. In such cases, an augmented Lagrange multiplier approach can provide a viable implementation (\cite{bertsekas2014constrained}).

\subsection{Structural Analysis} \label{subsec:XFEM}
% -------------------------------------------------------------------
The physical responses of the systems are predicted by the XFEM in this paper. A generalized Heaviside enrichment strategy (\cite{makhija2014numerical,terada2003finite,tran2011multiple}) is employed to avoid artificial coupling of disconnected material domains. The response is consistently interpolated in elemental subdomains with the same phase. The reader is referred to \cite{Barrera2020Topo} for more details on this topic. An in-depth description of the governing equations, including weakly enforcement of boundary conditions, and face-oriented ghost stabilization can be found in \cite{Barrera2019Hole}.
Since LS TO is not restricted to a particular immersed boundary technique, alternative approaches such as CutFEM (\cite{burman2015cutfem,burman2019cut}) or HIFEM (\cite{soghrati2016application}) could also be used for physical and sensitivity analyses.

% -------------------------------------------------------------------
% -------------------------------------------------------------------
\section{Numerical Examples}\label{sec:NumEx}
% -------------------------------------------------------------------
% -------------------------------------------------------------------

The proposed approach is studied using single material, solid-void linear elastic problems in 2D and 3D. Algorithmic parameter dependencies and performance are investigated with a typical structural beam problem. Strain energy minimization with mass constraint, and mass minimization with target strain energy optimization problem formulations are considered. The analyses provided in this section are focused on the minimum feature size control capability of the proposed approach. For details on algorithmic aspects of the hole nucleation strategy, the reader is referred to \cite{Barrera2019}.

The optimization problems are solved by the Globally Convergent Method of Moving Asymptotes (GCMMA, \cite{svanberg2002class}), with no inner iterations. The sensitivity analysis is performed via the adjoint method as detailed in \cite{sharma2017shape}. The optimization problem terminates when the constraints are satisfied and the relative change in the objective between three consecutive design iterations is less or equal to $1\times10^{-3}$. In all problems, the upper and lower bound of the LSF are function of the element size, $h$, and set to $\phi_{up} = 3h$ and $\phi_{low} = -3h$, respectively, unless specified otherwise. In all initial designs, a uniform LSF of $\phi(\mathbf{X})=\phi_0= 0.10\phi_{up}$ is prescribed, i.e., holes are not present a priori, unless specified otherwise. In addition, $\phi_0=\phi_{up}$, is prescribed in the vicinity of the Neumann and Dirichlet conditions to exclude them from the design domain. All simulations start with a uniform density of $\rho_{0}$, except the portions excluded from the design domain, which are prescribed a density of 1.0. Thus, a homogeneous porous material fills the entire design domain at the beginning of the optimization process. 

The continuation schemes described in Sections \ref{subsubsec:SIMPCont} and \ref{subsubsec:contRhoTh} are used to update the projection sharpness, $\gamma_{pr}$, the SIMP exponent, $\beta_{\rho}$, and the coupling density threshold, $\phi_{th_{hs}}$ using Eqs.~\ref{eq:projContFunc}, \ref{eq:simpExpContFunc}, and \ref{eq:expFuncContDec}, respectively. The update is performed every 10 optimization iterations in a total of 100 continuation steps. 
In addition, the VDDR penalty threshold is set to $\rho_{th_{fs}}=0.5\phi_{up}$, and the LS filter radius in Eq.~\ref{eq:LSLinearFilter1} is kept constant at $r_f^\phi$=1.5h and $r_f^\phi$=1.8h in 2D and 3D, respectively. Additional relevant parameters for the continuation schemes can be found in Table~\ref{tab:Ex1DefContParams}. Note that this continuation strategy was adopted for simplicity; however, more sophisticated schemes that update the algorithmic parameters based on the evolution of the design can be used.

% Table
\begin{table}[h]
	\caption{\label{tab:Ex1DefContParams}Default continuation and LS filter parameters for numerical examples.}
	\centering
	\renewcommand{\arraystretch}{1.2}
	\begin{tabular}{l|c}
		\hline
		Parameter                  & Value\\\hline
		Continuation step size&  $\mathcal{D}_{st}=10$              \\
		Number of design iterations in continuation &  $\mathcal{D}_{c}=100$              \\		
		Maximum number of design iterations &  $\mathcal{D}_{max}=175$              \\				
		Initial density threshold &  $\rho^0_{th}=0.25 \rho_0$              \\
		Initial SIMP exponent  &  $\beta_\rho^0=2.0$              \\
		Initial projection threshold&  $\gamma_{pr}^0=0.001$              \\	
		Final density threshold &  $\rho^f_{th_{hs}}=0.9 \rho_0$              \\
		Final SIMP exponent  &  $\beta_\rho^f=12.0$              \\
		Final projection threshold&  $\gamma_{pr}^f=40.0$              \\					
		Continuation density threshold exponent &  $\eta_{\rho_{th_{hs}}}=2.0$              \\
		Continuation SIMP exponent exponent &  $\eta_{\beta_{\rho}}=2.0$              \\		
		Continuation projection threshold exponent &  $\eta_{\gamma_{pr}}=2.0$       \\	
		\hline
	\end{tabular}
\end{table}
The XFEM is used to discretize the governing equations. The bulk material properties are shown in Table \ref{tab:LsDensCombProp} in self-consistent units. The physical and sensitivity analyses are performed using a domain decomposition approach and parallel computing. The linearized systems of equations are solved using the Multifrontal Massively Parallel Solver (MUMPS, \cite{amestoy2006hybrid}).
% Table
\begin{table}[h!]
	\caption{\label{tab:LsDensCombProp}Material properties interpolation parameters for numerical examples 1 and 2.}
	\centering
	\renewcommand{\arraystretch}{1.2}
	\begin{tabular}{l|c}
		\hline
		Property                  & Value\\\hline		
		Young's Modulus ($\Omega_{I}$)  	           &  $E_S = 2.0$x$10^3$              \\
		Young's Modulus ($\Omega_{II}$)                &  $E_V = 1$x$10^{-8}$ \\
		Poisson Ratio ($\Omega_{I}$ and $\Omega_{II}$) &  $\nu_S = \nu_V$ = 0.4 \\			
		Material Density ($\Omega_{I}$)     	           &  $\theta_S = 1.0$ \\
		Material Density ($\Omega_{II}$)              	   &  $\theta_V = 0.0$ \\
		\hline
	\end{tabular}
\end{table}
% Figure
\begin{figure}[h]
	\centering
	\includegraphics[width=0.5\linewidth]{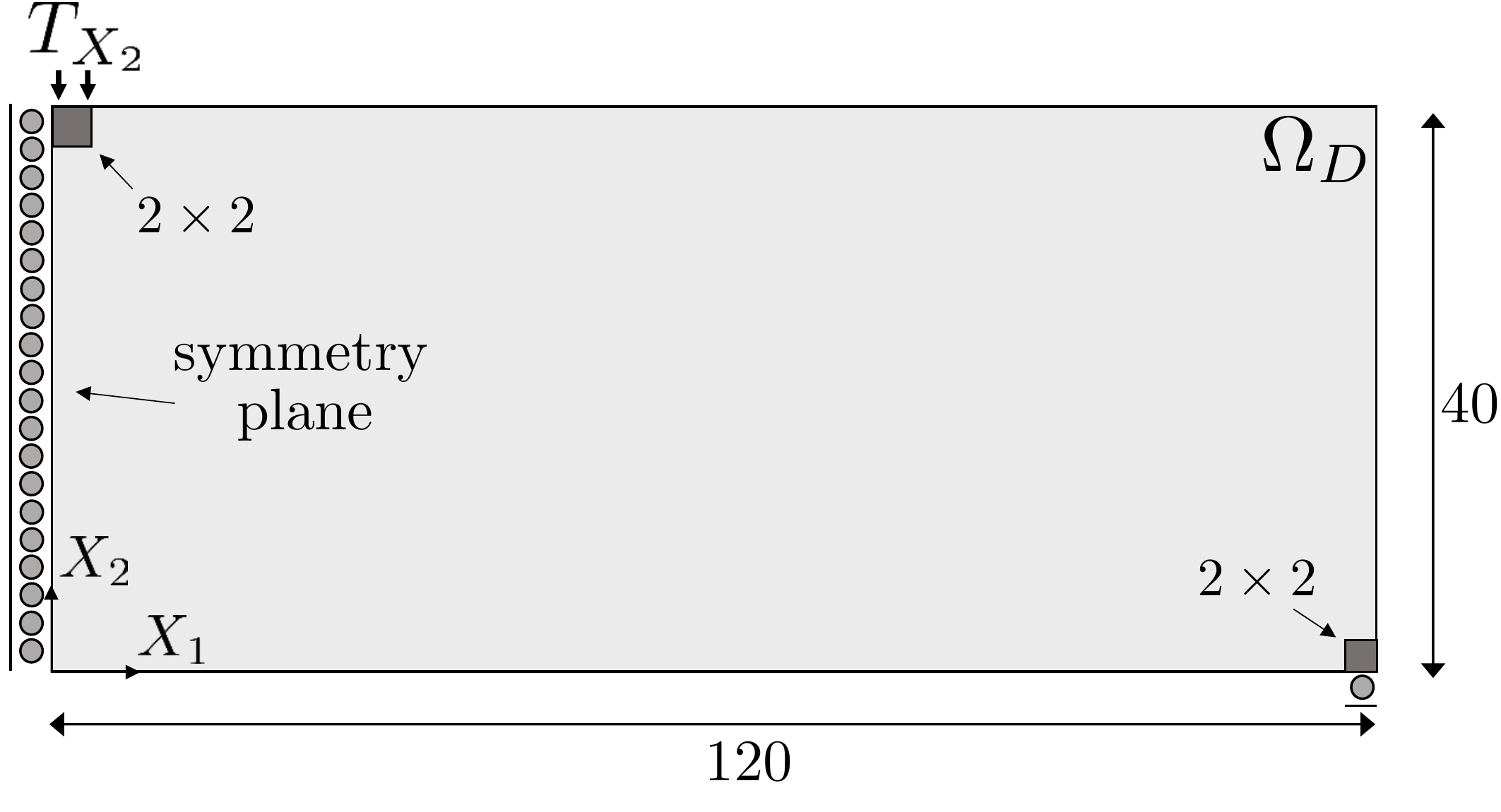}
	\caption{Problem setup of the 2D beam design problem. Only half of the domain is simulated considering symmetry.}
	\label{fig:beam2DProbSetup}
\end{figure} 

\subsection{Beam optimization (2D)} \label{sec:SupportStruc2D}
% -------------------------------------------------------------------------------
A beam optimization problem in 2D is studied to demonstrate the feature size control capability of the proposed framework. This first example focuses on examining the influence of algorithmic parameters of the projection and penalization of the density method, as well as the VDDR penalty. 
% Figure
\begin{figure*}[h!]
	\centering
	\includegraphics[width=1.0\linewidth]{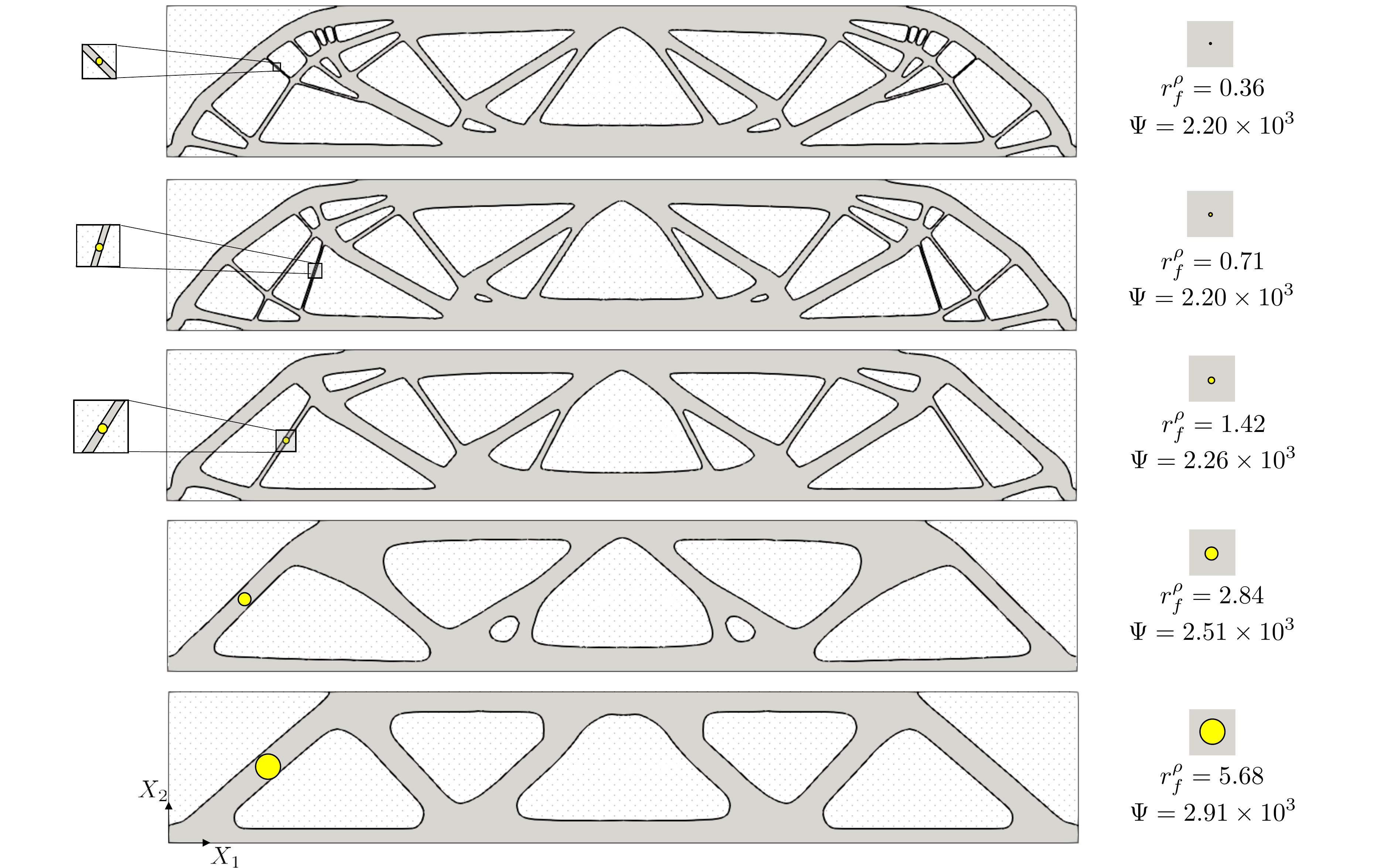}
	\caption{Optimized designs of 2D beam problem using different density filter radii.}
	\label{fig:finalDesMultFeatureSizes}
\end{figure*}
The loads and boundary conditions applied to this problems are specified in Fig.~\ref{fig:beam2DProbSetup}. A vertical traction of $T_{X_2} = -10.0$ is applied at the center of the top boundary of the design domain. Zero displacements are applied to the bottom corners in $X_2$ direction, and along the vertical axis (at $X_1=0$) in the $X_1$ direction of the design domain, to impose symmetry about $X_2$. Regions excluded from the design domain are highlighted in dark grey in Fig.~ \ref{fig:beam2DProbSetup}. Half of the domain of size $240\times40$ is discretized by a structured mesh with a uniform element size of $h=0.5$.
The strain energy minimization problem with mass constraint is formulated as follows:
% Equation
\begin{equation}\label{eq:StrEnMinProbSetup}
\begin{aligned}
\underset{s}{\min}~z(\mathbf{s}, \mathbf{u (\mathbf{s})})  
& = 
 w_1 ~ \frac{\Psi(\mathbf{s}, \mathbf{u(\mathbf{s})})}{\Psi_0}
+ 
w_2  ~ P_{Per}(\mathbf{s})
\\
& 
+ 
w_3 ~ P_{Reg}(\mathbf{s})
+ 
w_4 ~ P_{\hat\rho\phi}(\mathbf{s})
\\
& 
+
w_5 ~ P_{\hat\rho^{\Omega_{II}}}(\mathbf{s})
\\
s.t.: ~~~
&g_1(\mathbf{s})  = \frac{\mathcal{M}^{\Omega_1}(\mathbf{s})}{\mathcal{M}_{0}} - \gamma_m \leq 0.
\end{aligned}
\end{equation}

The strain energy, $\Psi$, is normalized by the strain energy of the initial design, $\Psi_0$. Details of the remaining objective components can be found in Section \ref{sec:LsOptFramework}. The initial objective weights are $w_i$ =[ 0.8345, 0.005, 0.015, 0.150, 1.50]. The perimeter penalty, $w_2$, is updated at every continuation step using the formulation in Eq.~\ref{eq:expFuncContDec} until it reaches 0.01. 
In constraint $g_1$, the mass of the solid phase, $\mathcal{M}^{\Omega_1}$, is normalized by a mass, $\mathcal{M}_{0}$, defined as the volume of the design domain, $\Omega_{D}$, multiplied by a density of 1.0 for unit consistency. The target mass is set to $\gamma_m = 0.40$.
% Figure
\begin{figure}[h]
	\centering
	\includegraphics[width=0.5\linewidth]{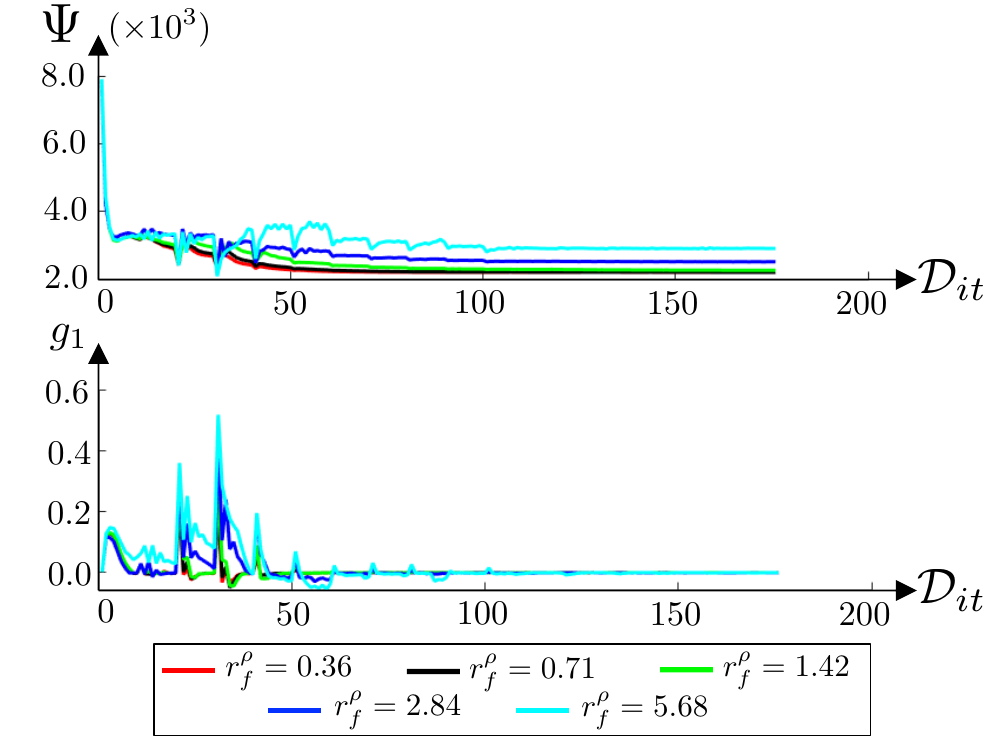}
	\caption{Strain energy and constraint evolution of 2D beam problem.}
	\label{fig:beam2DProbObjAndConstEvol}
\end{figure}
% Figure
\begin{figure*}[h!]
	\centering
	\includegraphics[width=1.0\linewidth]{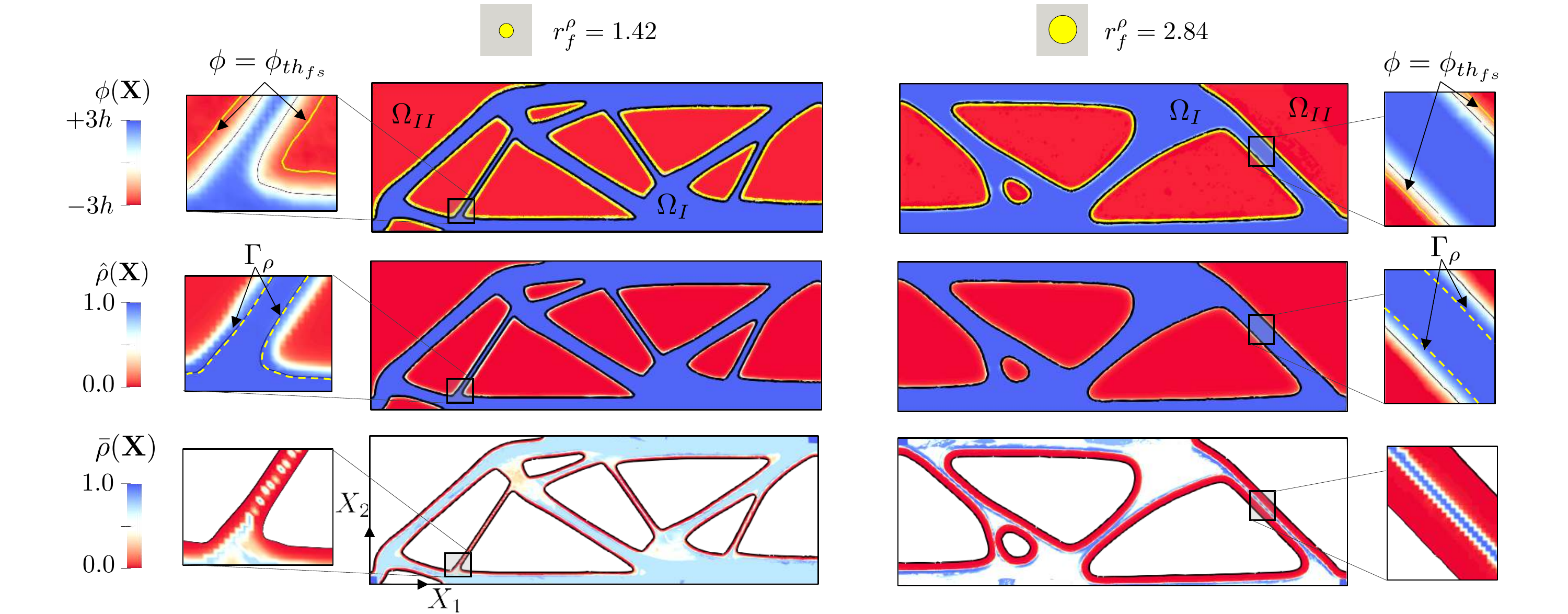}
	\caption{Half of the optimized designs of the the 2D beam optimization problem showing the (top) LS, (center) filtered projected density, and (bottom) unfiltered unprojected density fields for two different filter radii.}
	\label{fig:InvolvedDensityFields}
\end{figure*}

\subsubsection{Minimum feature size control} \label{subsec:2dBeamFSCResults}
% ------------------------------------------------------------------------------------------------------------
Figure~\ref{fig:finalDesMultFeatureSizes} shows the optimized designs of the 2D beam problem using five different density filter radii: $r_f^\rho$=[0.36, 0.71, 1.42, 2.84, 5.68]. 
The density field provides appropriate guesses for hole nucleation during the first stages of the optimization process. In addition, connections thinner than the density filter radius, $r_f^\rho$, are effectively removed in all optimized designs. The desired feature size control is achieved without changing the LS filter radius, $r_f^\phi$. As expected, there is a trade-off between increasing the minimum length allowed and the performance of the converged design. The larger $r_f^\rho$, the larger the final strain energy is.

The evolution of $\Psi$ and $g_1$ for all the problem setups of Fig.~\ref{fig:finalDesMultFeatureSizes} are depicted at the top and bottom of Fig.~\ref{fig:beam2DProbObjAndConstEvol}, respectively.
Considerable variations in both $\Psi$ and $g_1$ at the initial stages of the optimization process are observed at every update of the continuation scheme, where most of the holes are nucleated. Note that the effect of these jumps diminishes as the optimization problem progresses due to having a more uniformly defined density field in the solid phase. This is a direct result of the gradual increase of $\beta_\rho$ and $\gamma_{pr}$. A smooth evolution of both the objective and the mass constraint is observed at a later stage of the optimization process for all cases.

\subsubsection{Interplay of design fields} \label{subsec:2dBeamInterplay}
% ------------------------------------------------------------------------------------------------------------
Figure~\ref{fig:InvolvedDensityFields} shows the optimized designs of Fig.~\ref{fig:finalDesMultFeatureSizes} corresponding to $r_f^\rho$=[1.42, 2.84] colored by three fields: the LSF, projected density field, and an unfiltered unprojected density field, $\bar\rho(\mathbf{X})$. This last field is constructed by linearly interpolating the density optimization variables as follows:
% Equation
\begin{equation}\label{eq:discrLsField2}
\begin{aligned}
\bar\rho(\mathbf{X})
= 
\displaystyle\sum_{i=1}^{N^e} \mathcal{N}_i ( \mathbf{X}) ~ s_i^\rho.
\end{aligned}
\end{equation}
In addition, insets with the contours described in Section \ref{subsec:fscStrategy} are shown to assess the effectiveness of the minimum feature size functionality.

The top row of Fig.~\ref{fig:InvolvedDensityFields} shows that the LSF has a smooth transition between material and void with a unit gradient norm, as well as a uniform truncated field between [+3h,-3h] away from the interface.  This is achieved by the regularization penalty summarized in Section \ref{sec:LsOptFramework}. The yellow contour $\phi(\mathbf{X})=\phi_{th_{fs}}=1.5h$ follows closely the interface while leaving room for intermediate densities to develop in the space between this contour and the interface. A smooth LSF is observed with both filter radii. As explained in Section \ref{subsec:VDDRPen}, failing to enforce a regularized LSF like the one shown in this example prevents the control on the region where the VDDR penalty is active in the void phase.
% Figure
\begin{figure}[h]
	\centering
	\includegraphics[width=0.5\linewidth]{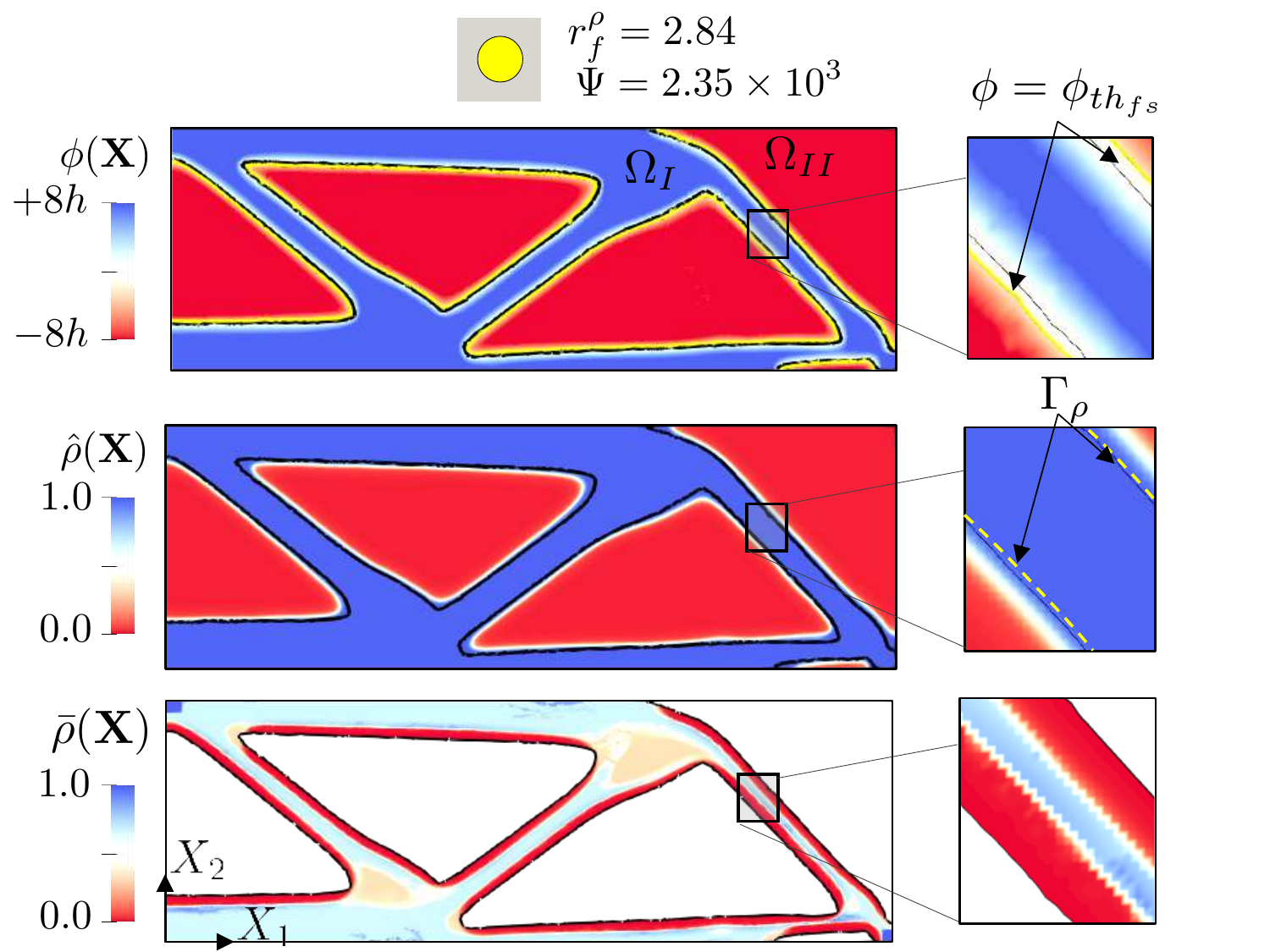}
	\caption{Half of the optimized designs of the 2D beam optimization problem showing the (top) LS, (center) filtered projected density, and (bottom) unfiltered unprojected density fields using $r_f^\rho$=2.84, $\phi_{up}=8.0h$ and $\phi_{low}=-8.0h$.}
	\label{fig:beam2DInfluenceLSRegBonds}
\end{figure} 

The projected density fields of the optimized designs shown in the middle row of Fig.~\ref{fig:InvolvedDensityFields} illustrate fields predominantly filled with uniform densities of 0.0 and 1.0 in the void and solid phases, respectively. Insets of this field display narrow transition regions for both density filter radii. It can be seen that for $r_f^\rho$=1.42 the interface and the density contour $\Gamma_\rho$ (depicted using yellow dashed lines) overlap, and thus achieve the desired behavior; see Fig.~\ref{fig:DensityIssueInVoidRegion}. However, overcutting, i.e., placing low-density material within the solid phase in the vicinity of the interface is observed for $r_f^\rho$=2.84. In general, the larger $r_f^\rho$ is, the more diffuse is the transition region of the filtered projected density field. Despite the projection promoting a narrow transition region, its width strongly depends on the filter radius. Note that overcutting can be mitigated by choosing the bounds of the truncated regularized LSF and the $\phi_{th_{fs}}$ threshold as function of $r_f^\rho$ appropriately. This way, material is not removed too closely to the interface. Figure~\ref{fig:beam2DInfluenceLSRegBonds} shows this effect on the results with $r_f^\rho$=2.84 by using $\phi_{up}=8h$ and $\phi_{low}=-8h$ as the bounds of the regularized LSF. 
% Figure
\begin{figure*}[h!]
	\centering
	\includegraphics[width=1.0\linewidth]{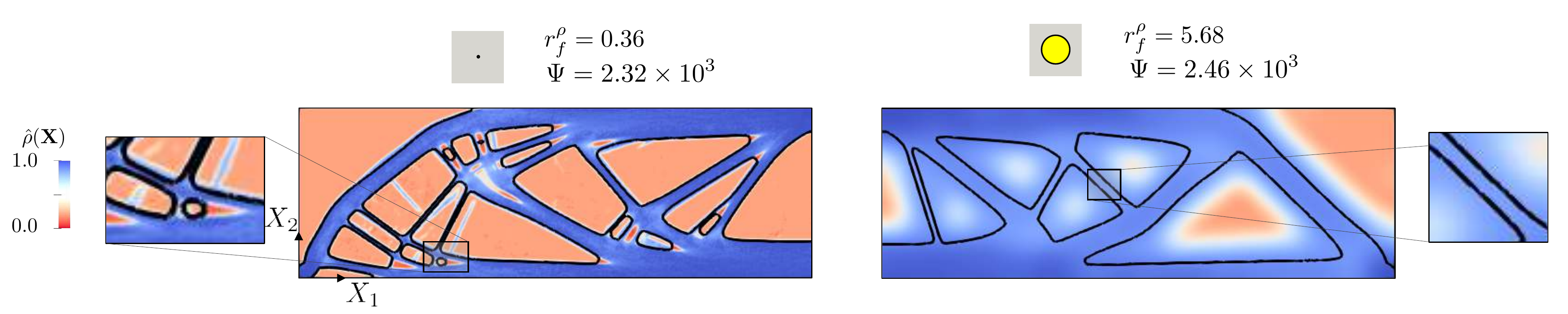}
	\caption{Optimized design of 2D beam problems using two different filter radii and a setup without the components that provide minimum feature size control.}
	\label{fig:finalDesNoFeatureSizeControl}
\end{figure*}

The unfiltered unprojected density field, $\bar\rho(\mathbf{X})$, is shown at the bottom row of Fig.~\ref{fig:InvolvedDensityFields}. Only the solid phase is displayed for clarity but a uniform field, $\bar\rho(\mathbf{X})\approx0.0$, is obtained in the void phase. The insets show a skeleton delineated at the center of regions where the minimum feature size is enforced in both cases. This is achieved by promoting a density of zero in the void phase through the VDDR penalty. Note that the skeleton is better resolved for the results with $r_f^\rho$=2.84 due to the VDDR penalty pushing the transition region to the solid phase. However, if the density field on the solid phase is uniform, as seen in Figure~\ref{fig:beam2DInfluenceLSRegBonds}, the skeleton is characterized by intermediate unfiltered unprojected densities.

\subsubsection{Comparison against results without minimum feature size control} \label{subsec:2dBeamComparison}
% ------------------------------------------------------------------------------------------------------------
The goal of this section is to show that the components considered in the proposed TO approach, i.e., large SIMP exponent and projection sharpness, as well as a VDDR penalty are indeed essential to achieve minimum feature size control. To demonstrate this, a similar problem setup with a (low) constant SIMP exponent ($\beta_\rho^0=\beta_\rho^f$=2.0) no projection ( $\gamma_{pr}^0=\gamma_{pr}^f$=0.0001), and no VDDR penalty ($w_5=0.0$) is considered. Both the smallest and largest density filter radii of the configurations shown in Fig~\ref{fig:finalDesMultFeatureSizes}, namely $r_f^\rho$=[0.36, 5.68], are used.

The optimized designs are shown in Fig.~\ref{fig:finalDesNoFeatureSizeControl}. It can be seen that increasing $r_f^\rho$ has a noticeable effect on the final topologies. This is a consequence of the employed hole nucleation strategy; see Section \ref{subsec:holeSeeding}. A larger $r_f^\rho$ reduces the number of holes seeded through the optimization process. Thus, the problems converge to different local optima. Nonuniform density distributions in both phases and a lack of control of the minimum features of the optimized designs are observed in both cases. On the left side of Fig~\ref{fig:finalDesNoFeatureSizeControl}, the solid phase of the optimized design is characterized by multiple regions with intermediate and low densities. Furthermore, the projected density field in the void phase settles to $\hat\rho\approx0.30$ since no VDDR penalty is enforced. On the right side of Fig~\ref{fig:finalDesNoFeatureSizeControl}, considerable undercutting by the LSF is observed. In addition, the design domain is predominantly populated by intermediate densities. In both cases, a simplified problem setup that relies solely on the density filter radius fails to provide explicit minimum feature size control. 

\subsection{Beam optimization (3D)} \label{sec:SupportStruc3D}
\begin{figure}[h]
	\centering
	\includegraphics[width=0.5\linewidth]{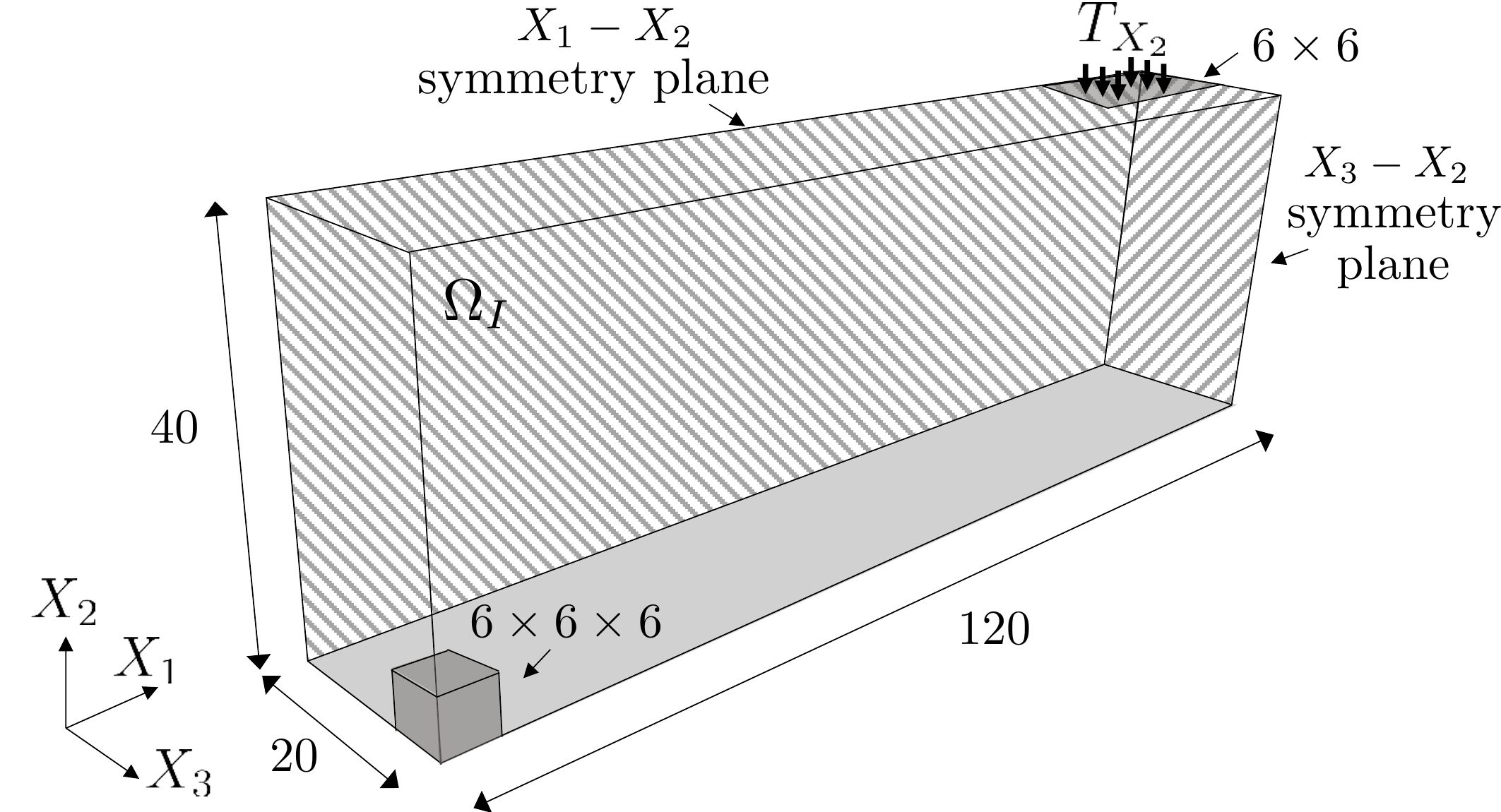}
	\caption{Quarter of design domain with boundary conditions, symmetry planes and dimensions.}
	\label{fig:3dMbbBeamProbSetup}
\end{figure}

% -------------------------------------------------------------------------------
An equivalent 3D setup of the 2D beam described above is studied next. The behavior of the proposed approach is examined using two optimization problem formulations. Mesh independence of the proposed approach is studied using a strain energy minimization with mass constraint formulation. In addition, an alternative to the VDDR penalty to promote low densities in the void phase is investigated in a mass minimization with strain energy constraint optimization problem.
A design domain of $240\times40\times40$ is simply supported at all corners of the bottom face. In addition, a traction load $T_{X_2} = -10.0$ is applied on a rectangular region at the center of the top face. Only one quarter of the design domain is analyzed and optimized to exploit symmetry. The 3D beam problem characteristics are shown in Fig. \ref{fig:3dMbbBeamProbSetup}. The quarter domain of size $120\times40\times20$ is discretized by a structured mesh. The material properties assigned to both phases are the same as in the previous example; see Table \ref{tab:LsDensCombProp}.
% Figure
\subsubsection{Strain energy minimization: hole seeding and minimum feature size control through density method} \label{subsec:3dBeamHSnFSC}

\begin{figure*}[h!]
	\centering
	\includegraphics[width=0.99\linewidth]{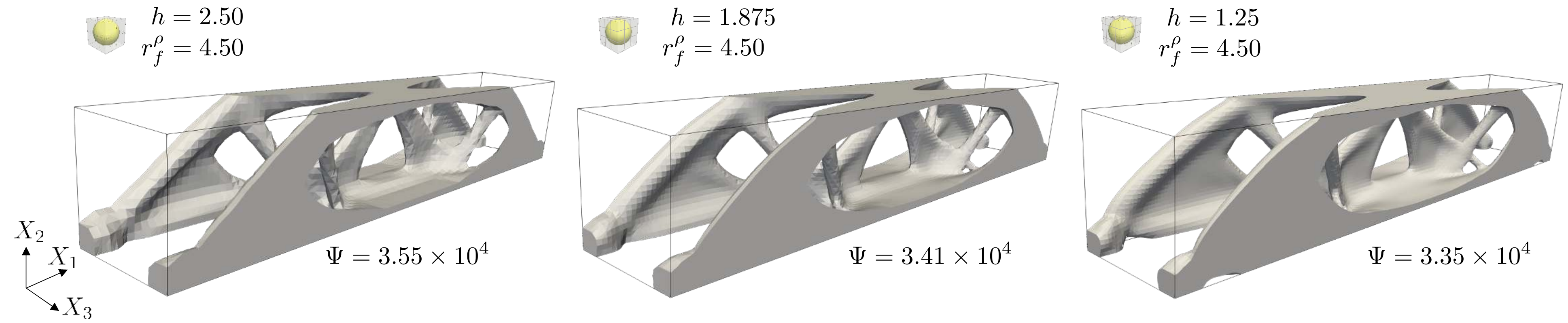}
	\caption{Optimized 3D beams using three mesh sizes for a strain energy minimization with mass constraint problem.}
	\label{fig:finalDesStrainEnProbMeshRef}
\end{figure*}

Mesh-independence of the proposed approach is studied in this section using three meshes with element edge lengths of $h=[2.50, 1.875, 1.25]$.
The optimization problem formulation of Eq.~\ref{eq:StrEnMinProbSetup} is used. The objective weights in this setup are $w_i$ =[ 0.739, 0.001, 0.010, 0.250, 2.500]. A mass constraint of $\gamma_m = 0.15$ is enforced.
The optimized structures for all three levels of refinement are summarized in Fig.~\ref{fig:finalDesStrainEnProbMeshRef}. The same density filter radius, i.e., $r^\rho_f=4.50$, is used in all simulations.
Similar performances are achieved, and the converged topologies differ only locally. Truss-like structures with partial shear webs are obtained.
These results suggest that the TO approach introduced in this paper leads to mesh-independent converged designs, at least for the class of optimization problems studied by the numerical examples studied. 
% Figure

\subsubsection{Strain energy minimization: minimum feature size control through density method without hole seeding} \label{subsec:3dBeamHSnFSC2}
\begin{figure*}[h!]
	\centering
	\includegraphics[width=0.99\linewidth]{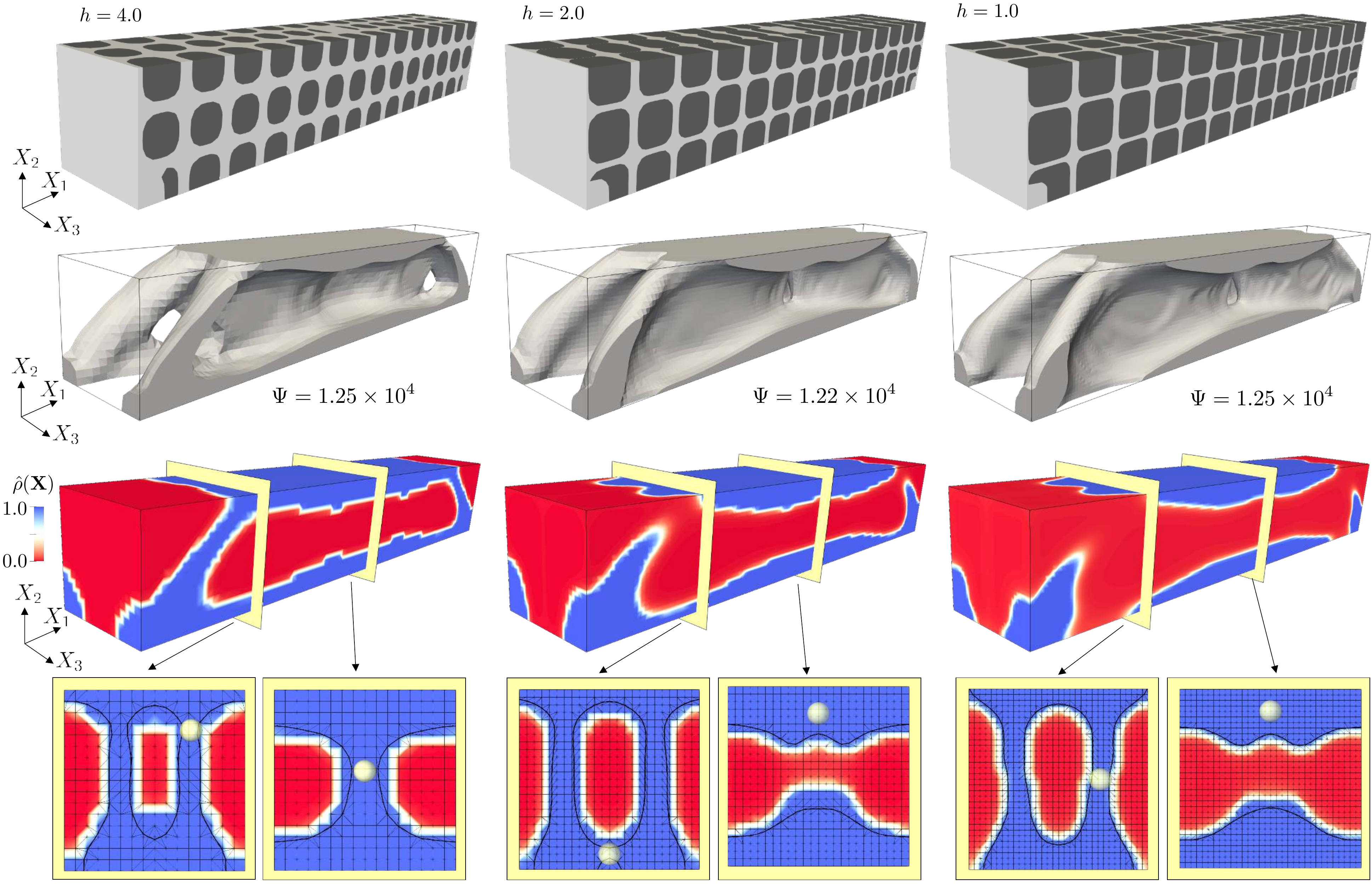}
	\caption{3D beam strain energy minimization with mass constraint problem using three mesh sizes and initial hole pattern: (top) initial design domains, (center) solid phase of converged designs, and (bottom) design domain at convergence colored by projected density field with insets.}
	\label{fig:finalDesStrainEnProbMeshRef_NoHoleSeedPen}
\end{figure*}

% ------------------------------------------------------------------------------------------------------------
To demonstrate that the minimum length control is completely decoupled from the hole seeding strategy, a modified version of the optimization formulation of Eq.~\ref{eq:StrEnMinProbSetup} is investigated. The same meshes of the previous setup are used, i.e. $h=[2.50, 1.875, 1.25]$.
In this case, the hole seeding penalty term is deactivated by setting $w_4$ to zero. Instead, the initial designs contain a hole seeding pattern of cuboid primitives, as seen in the first row of Fig.~\ref{fig:finalDesStrainEnProbMeshRef_NoHoleSeedPen}.
In this setup, the projected density field is initially set to 1.0 in the entire design domain and can freely evolve throughout the optimization process.

The second row of Fig.~\ref{fig:finalDesStrainEnProbMeshRef_NoHoleSeedPen} shows the optimized designs. The final topologies consist of the typical shear-webs observed in 3D benchmark LS TO beam problems. Truss-like connections are not observed since holes are no longer seeded through the optimization process, and the projected density field is used only to provide minimum feature size control. Overall, the optimized structures differ only locally, and similar minimum length is achieved. The difference in converged designs is less pronounced between the two finest meshes since the difference in the initial hole seeing patterns is small. Note that identical initial hole seeding patterns cannot be achieved because the hole arrangement is constraint by the element length, $h$. Placing holes closer than done in the first row of Fig.~\ref{fig:finalDesStrainEnProbMeshRef_NoHoleSeedPen} may result in spurious connections between them; see, for example, \cite{Barrera2020Topo}.

The design domain colored by the projected density field is shown at the bottom row of Fig.~\ref{fig:finalDesStrainEnProbMeshRef_NoHoleSeedPen} together with two views of 2D cuts in the $X_2-X_3$ plane. These cuts are located at $X_1=75.0$ and $X_1=150.0$ and show the mesh along with a yellow sphere which radius is equal to the prescribed density filter radius. The spheres are placed where the minimum features are located. It can be seen that the VDDR penalty effectively enforces a density of zero in the void phase, and neither undercutting nor overcutting is observed for any of the refinement levels. In addition, a uniform projected density field of 1.0 is observed in the solid phase. As expected, the resolution of the transition region increases as the mesh is refined and, in all cases, it is located close to the interface in the void phase, achieving the desired configuration shown in Fig.~\ref{fig:DensityIssueInVoidRegion}. In general, a consistent minimum length scale control is achieved by the proposed TO approach regardless of whether the projected density field is used for hole nucleation or not.
% Figure

\subsubsection{Total mass minimization} \label{subsec:massMin3DBeam}
% ------------------------------------------------------------------------------------------------------------
An optimization problem formulation which incorporates the functionality of the VDDR penalty in the void phase is investigated in this section. 
The proposed mass minimization problem with a strain energy constraint reads:
% Equation
\begin{equation}\label{eq:MassMinProbSetup}
\begin{aligned}
\underset{s}{\min}~z(\mathbf{s}, \mathbf{u (\mathbf{s})})  
& = 
 w_1 ~  \frac{\mathcal{M(\mathbf{s})}}{\mathcal{M}_{0}} 
+ 
w_2  ~ P_{Per}(\mathbf{s})
\\
& 
+ 
w_3 ~ P_{Reg}(\mathbf{s})
+ 
w_4 ~ P_{\rho\phi}(\mathbf{s})
\\
&
+
w_5 ~ P_{\hat\rho^{\Omega_{II}}}(\mathbf{s})
\\
s.t.: ~~~
&g_1(\mathbf{u}, \mathbf{s})  = \frac{\Psi(\mathbf{s}, \mathbf{u(\mathbf{s})})}{\Psi_{ref}} - 1.0 \leq 0.
\end{aligned}
\end{equation}
The first term in the objective is the mass of the entire design domain, $\mathcal{M}$, defined as:
% Equation
\begin{equation}\label{eq:massMinObj}
\begin{aligned}
\mathcal{M}(\mathbf{s})
= 
\displaystyle\int_{\Omega_{I}} \Theta dV
+
\displaystyle\int_{\Omega_{II}} \Theta dV,
\end{aligned}
\end{equation}
which indicates that the mass is minimized in both the solid and void phases. Note that the mass constraint in the previous optimization problem only considers the solid phase, see Eq.~\ref{eq:StrEnMinProbSetup}. The remaining components are explained in Section \ref{sec:LsOptFramework}. The initial objective weights are $w_i$ =[0.819, 0.001, 0.03, 0.15, 0.00]. Note that $w_5$ is set to zero to deactivate the VDDR penalty. A strain energy constraint that is normalized by $\Psi_{ref}=1.5\times 10^5$ is enforced.

\begin{figure*}[t]
	\centering
	\includegraphics[width=1.0\linewidth]{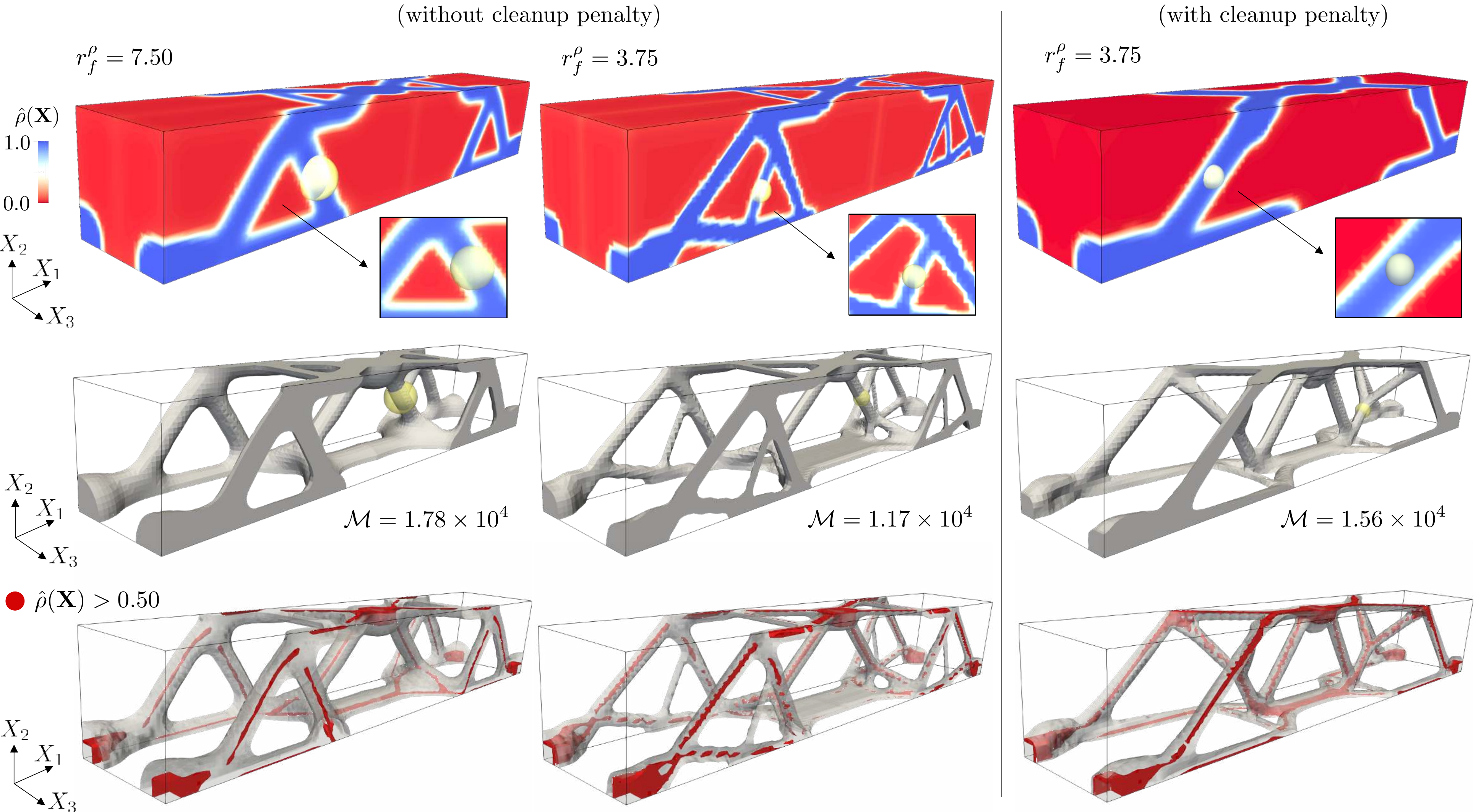}
	\caption{3D beam mass minimization problem. Showing: (top) only solid domain, (center) the entire design domain colored by the filtered projected density field, and (bottom) the solid domain highlighting in red unfiltered unprojected densities greater than a threshold.}
	\label{fig:finalDesSkelMultFiltRad}
\end{figure*}

The first, second and third row of Fig.~\ref{fig:finalDesSkelMultFiltRad} display the entire design domain colored by the projected density field, the optimized designs together with a yellow sphere that represents the desired minimum feature, and the skeleton of the converged designs in the solid phase, respectively. The results for this new setup are shown in the left and middle columns, for $r_f^\rho=[7.50,3.75]$. The first row shows that a uniform density of $0.0$ is achieved in the void phase without including a VDDR penalty in the optimization problem formulation. However, the explicit control of the effect of the VDDR penalty, i.e., how far from the interface it is enforced, is lost. This has a number of implications: (i) the transition region develops within the bounds of the filter radius; (ii) the interface enforces a minimum feature drastically different from the desired one, and its effect is more pronounced with larger filter radii (see second row of Fig.~\ref{fig:finalDesSkelMultFiltRad}); and (iii) a clearly defined skeleton is less likely to occur (see third row of Fig.~\ref{fig:finalDesSkelMultFiltRad}).
% Figure

To show the effect of having explicit control of the minimum feature size, the optimization problem formulated in Eq.~\ref{eq:MassMinProbSetup} is solved with a VDDR penalty weight different than zero, i.e. $w_5=10.0$. Also, the mass is defined on the solid side only as follows:
% Equation
\begin{equation}\label{eq:massMinObj2}
\begin{aligned}
%\mathcal{M}_{\Omega_I}(\mathbf{s})
\mathcal{M}(\mathbf{s})
= 
\displaystyle\int_{\Omega_{I}} \Theta dV
=
\mathcal{M}^{\Omega_1}(\mathbf{s}).
\end{aligned}
\end{equation}
The right column of Fig.~\ref{fig:finalDesSkelMultFiltRad} shows the optimized design of this modified problem using a filter radius of $r^\rho_f=3.75$. It can be seen that a reliable minimum feature size is recovered. Thus, the VDDR penalty is crucial to achieve the desired control on the minimum feature.

% -------------------------------------------------------------------
% -------------------------------------------------------------------
\section{Conclusions and Future Work} \label{sec:Concl}
% --------------------------------------------------------------------
% -------------------------------------------------------------------

A TO approach that provides mechanisms for continuous hole seeding and minimum feature size control was introduced. This approach combines LS and density methods. LS and density optimization variables are coupled directly through a penalty term in the objective function to nucleate holes at the initial stages of the optimization process. In addition, the projection and penalization of the density field are gradually increased through the optimization process, and the density field in the void phase is kept near zero through a VDDR penalty added to the objective function. The minimum feature size control is provided by the density filter radius. This has the advantage that thin members can be removed during the optimization process if physically meaningful and changes in topology are not prohibited, as seen by typical LS-based feature size control schemes. In the proposed approach, the LS and density fields are optimized simultaneously throughout the entire optimization process. Hence, this approach has the flavor of a combined LS-density TO method. Despite projecting the density field, converged designs without intermediate densities cannot be avoided. 

Numerical examples in 2D and 3D demonstrated the capabilities of the proposed approach. Strain energy and mass minimization optimization problems were considered. It was shown that the VDDR penalty, which is a function of a LS threshold, requires the LS field to be regularized (i.e., use a method to avoid locally too flat or steep LSF in the vicinity of the interface) to effectively control the minimum feature size. Otherwise, this penalty may not be activated sufficiently away from the interface in the void phase. Mesh-independent designs were obtained for typical structural beam problems with minimized strain energy and mass, with and without hole seeding capabilities. Due to the number of components and algorithmic parameters needed in this framework, its applicability requires an in-depth understanding of the design optimization problem at hand. 

Future work should study the proposed approach for setups with intricate optimization problem formulations and different physics models to assess complexities associated with achieving well-posed optimization problem formulations. Also, more sophisticated continuation schemes that update the algorithmic parameters based on the evolution of the design should be explored.

% -------------------------------------------------------------------
% -------------------------------------------------------------------
\section*{Acknowledgements}
% -------------------------------------------------------------------
% -------------------------------------------------------------------

All authors acknowledge the support of the National Science Foundation (CMMI-1463287). 
The first author acknowledges partial auspice of the U.S. Department of Energy by Lawrence Livermore National Laboratory under Contract DE-AC52-07NA27344 (LLNL-JRNL-820531).
The third author acknowledges the support of the Defense Advanced Research Projects Agency (DARPA) under the TRADES program (agreement HR0011-17-2-0022). The opinions and conclusions presented in this paper are those of the authors and do not necessarily reflect the views of the sponsoring organizations.


% \bibliographystyle{abbrvnat}   % 
% \bibliographystyle{plainnat}   % 
% \bibliography{Mendeley}   % name your BibTeX data base
% \bibliography{holeSeedingBiblio}   % name your BibTeX data base

\end{document}